\documentclass[12pt]{article} 
 \usepackage[margin=1in]{geometry}
\usepackage{amssymb,amsfonts,amscd,amsthm}
\usepackage{amsmath}
\usepackage[all,arc]{xy}
\usepackage{enumerate, bbm}
\usepackage{mathrsfs}
\usepackage{tikz}
\usetikzlibrary{positioning}
\usetikzlibrary{calc}

\usepackage{mathtools}
\usepackage{hyperref}
\usepackage{color}
\usepackage{graphicx}
\usepackage[shortlabels]{enumitem}
\usepackage{enumitem} 
\graphicspath{ {TexImages/} }
\usepackage{cleveref}
\usepackage{thm-restate}

\newtheorem{thm}{Theorem}[section]
\newtheorem{cor}[thm]{Corollary}
\newtheorem{prop}[thm]{Proposition}
\newtheorem{obs}[thm]{Observation}
\newtheorem{lem}[thm]{Lemma}

\newtheorem{defn}[thm]{Definition}
\newtheorem{claim}{Claim}

\newtheorem{const}{Construction}

\newcommand{\floor}[1]{\left\lfloor #1\right\rfloor}

\definecolor{amber}{rgb}{1.0, 0.75, 0.0}

\definecolor{forest}{rgb}{0.0, 0.5, 0.0}

\definecolor{cadmium}{rgb}{0.93, 0.53, 0.18}

\definecolor{bittersweet}{rgb}{1.0, 0.44, 0.37}

\definecolor{byzantine}{rgb}{0.74, 0.2, 0.64}

\definecolor{brilliantrose}{rgb}{1.0, 0.33, 0.64}

\DeclareMathOperator{\ex}{ex}
\DeclareMathOperator{\sat}{sat}
\DeclareMathOperator{\rsat}{sat^*}

\title{On the proper rainbow saturation numbers of cliques, paths, and odd cycles}
\author{Dustin Baker\thanks{Department of Mathematics, Iowa State University, Ames IA.} \and Enrique Gomez-Leos\footnotemark[1] \and Anastasia Halfpap\footnotemark[1] \and Emily Heath\thanks{Department of Mathematics and Statistics, California State Polytechnic University, Pomona, CA.} \and Ryan R. Martin\footnotemark[1] \and Joe Miller\footnotemark[1] \and Alex Parker\footnotemark[1] \and Hope Pungello\footnotemark[1] \and Coy Schwieder\footnotemark[1] \and Nick Veldt\footnotemark[1]}
\date{\today}

\begin{document}

\maketitle

\begin{abstract}
Given a graph $H$, we say a graph $G$ is \emph{properly rainbow $H$-saturated} if there is a proper edge-coloring of $G$ which contains no rainbow copy of $H$, but adding any edge to $G$ makes such an edge-coloring impossible. The \emph{proper rainbow saturation number}, denoted $\rsat(n,H)$, is the minimum number of edges in an $n$-vertex rainbow $H$-saturated graph. We determine the proper rainbow saturation number for paths up to an additive constant and asymptotically determine $\rsat(n,K_4)$. In addition, we bound $\rsat(n,H)$ when $H$ is a larger clique, tree of diameter at least 4, or odd cycle. 
\end{abstract}

\section{Introduction}

For a fixed graph $H$, how many edges can a graph $G$ on $n$ vertices have if it does not contain $H$ as a subgraph? We say that $G$ is \emph{$H$-saturated} if $G$ contains no copy of $H$, but for any $x,y\in V(G)$ with $xy\notin E(G)$, the graph $G+xy$ on vertex set $V(G)$ with edge set $E(G)\cup\{xy\}$ contains a copy of $H$. 

A classical question in extremal combinatorics asks for the maximum number of edges in an $n$-vertex $H$-saturated graph. This is called the \emph{Tur\'an number} $\ex(n,H)$, and has been extensively studied following work of Mantel \cite{Mantel} and Tur\'an \cite{Turan41} determining $\ex(n,K_t)$ for $t\geq 3$. Also of interest is the other extremal case, in which we seek the minimum number of edges in an $n$-vertex $H$-saturated graph, called the \emph{saturation number} $\sat(n,H)$. The study of the saturation number was initiated by work of Zykov \cite{zykov} and independently Erd\H{o}s, Hajnal, and Moon \cite{erdos1964problem}. Many different generalizations of these two problems have been studied over the years, including analogous questions in the setting of edge-colored graphs.

An edge-coloring $c:E(G)\rightarrow C$ is \emph{proper} if $c(e)\neq c(f)$ for all incident edges $e,f$ and \emph{rainbow} if $c(e)\neq c(f)$ for all edges $e,f\in E(G)$. Given graphs $G$ and $H$ and a proper edge-coloring $c$ of $G$, we say that $G$ is \emph{rainbow $H$-free under $c$} if $G$ contains no copy of $H$ which is rainbow under the coloring $c$. A graph $G$ is \emph{(properly) rainbow $H$-saturated} if the following two conditions hold:

\begin{enumerate}
    \item There exists a proper edge-coloring $c$ of $G$ such that $G$ is rainbow $H$-free under $c$, and
    \item For any edge $e\notin E(G)$, any proper edge-coloring of $G+e$ contains a rainbow copy of $H$.
\end{enumerate}

Keevash, Mubayi, Sudakov, and Verstra\"ete \cite{KMSV} introduced the \emph{rainbow extremal number} $\ex^*(n,H)$ which is the maximum number of edges in an $n$-vertex rainbow $H$-saturated graph. More recently,  Bushaw, Johnston, and Rombach \cite{Bushaw2022} initiated the study of the \emph{proper rainbow saturation number} $\rsat(n,H)$ which is the minimum number of edges in an $n$-vertex rainbow $H$-saturated graph. 

It is worth noting that other versions of rainbow saturation problems have been considered as well. For example, Behague, Johnston, Letzter, Morrison, and Ogden \cite{behague2022rainbow} explored the problem without the restriction that colorings be proper. Another variant, introduced by Hanson and Toft \cite{hanson1987}, requires that colorings, while not necessarily proper, are restricted to a set of only $t$ colors. 

Throughout the remainder of the paper, we will exclusively study proper rainbow saturation problems; thus we drop the qualifier ``proper" as this is assumed. 

There are few graphs for which the rainbow saturation number is known asymptotically, namely $\rsat(n,P_4)$ determined by Bushaw, Johnston, and Rombach \cite{Bushaw2022} and $\rsat(n,C_4)$  determined by Halfpap, Lidick\'y, and Masa\v{r}\'ik \cite{halfpap2024proper}. Adding to this list, we determine $\rsat(n,K_4)$ asymptotically. 

\begin{restatable}{thm}{kfourlowerbound}\label{thm: k4 lower bound}
Let $\frac{1}{2} > \alpha > 0$. For any $n$ such that $\alpha^2 n \geq 7$ and $\alpha n > 220$, we have
\[\frac{7}{2}n - 8 \alpha n\leq \rsat(n,K_4)\leq \frac{7}{2}n+O(1).\]
\end{restatable}

We also give bounds on the rainbow saturation number for larger cliques. 

Our next main result is to determine the rainbow saturation numbers  for paths up to an additive constant. Throughout the paper, we use $P_k$ to denote the path on $k$ vertices (with $k-1$ edges). 

\begin{restatable}{thm}{pathsthm}\label{paths thm}
For $k\geq 5$ and $n \geq (k-1)2^{k-4}$, we have 
\[n-1\leq \rsat(n,P_k)\leq n + O(2^k).\]
\end{restatable}

In particular, taking $k$ fixed relative to $n$, we have that $\mathrm{sat^*}(n,P_k)$ is asymptotically equal to $n$. We also remark that the lower bound of Theorem~\ref{paths thm} can be extended to a larger class of graphs, showing that $\rsat(n,T)\geq n-1$ for any tree $T$ of diameter at least 4. 

Finally, we consider the rainbow saturation number of cycles. Recently, Halfpap, Lidick\'y, and Masa\v{r}\'ik \cite{halfpap2024proper} proved that $\rsat(n,C_4)=\frac{11}{6}n+o(n)$. They also gave upper bounds for other short cycles, showing that $\rsat(n, C_5) \leq \left\lfloor\frac{5n}{2}\right\rfloor-4$ for $n\geq 9$  and that $\rsat(n, C_6) \leq\frac{7}{3}n+O(1)$. 
In the following theorem, we give an upper bound for the rainbow saturation number for longer odd cycles. For context, no disconnected graph can be rainbow $C_k$-saturated for any $k$, so we trivially have $n-1 \leq \mathrm{sat}^*(n,C_k)$ for all $k$.

\begin{restatable}{thm}{cyclesthm}\label{cycles thm}    
For odd $k \geq 7$ and for $n \geq 3k - 2$, we have 
\[\sat^*(n,C_k) \leq \left( \frac{k - 1}{2} \right) n -\binom{\frac{k + 1}{2}}{2}.\] 
\end{restatable}

The remainder of the paper is organized as follows. In Section \ref{sec:cliques}, we prove Theorem~\ref{thm: k4 lower bound} and bounds for larger cliques. In Section \ref{sec:paths}, we prove Theorem~\ref{paths thm}. Finally, in Section \ref{sec:cycles}, we prove Theorem~\ref{cycles thm}. 

\vspace{.5cm}
\textbf{Remark.} While completing this paper, we learned that Lane and Morrison simultaneously and independently derived partially-overlapping results for the proper rainbow saturation number of various graphs. 
Among other results in \cite{MLgeneral}, they gave the upper bound (but not the lower bound) on $\rsat(n,K_4)$ in Theorem \ref{thm: k4 lower bound} as well as more general upper bounds for larger cliques. They also obtained the upper bound on $\rsat(n,C_k)$ for odd $k$ in Theorem \ref{cycles thm} and gave a similar upper bound for even $k$. In \cite{MLtrees}, they used a similar approach to ours to show that $\rsat(n,P_k)=n+O(k)$ and studied other families of trees such as brooms and caterpillars which we did not consider. 

\subsection{Notation and Preliminary Definitions}
Throughout the paper, we will use the following notation. 
We denote the \textit{degree} of a vertex $v \in V(G)$ by $d(v)$ and the minimum vertex degree of a graph $G$ by $\delta(G)$. For a vertex $v\in V(G)$ we use $N(v)$ to denote the \textit{neighborhood} of $v$ by $N(v):=\{u \in V(G) : uv\in E(G) \}$ so that $d(v) = |N(v)|$. We also use $N[v]$ to denote the \textit{closed neighborhood} of $v$, that is, $N[v]:= N(v) \cup \{ v\}$. Given graphs $G$ and $H$, let $G \vee H$ denote the \emph{join} of $G$ and $H$ which has vertex set $V(G)\cup V(H)$ and edge set $E(G)\cup E(H)\cup \{gh:g\in V(G), h\in V(H)\}$. Let $E_n$ denote the empty graph on $n$ vertices. Given two sets $A,B$ we denote by $A \Delta B$ the \textit{symmetric difference}, that is $A \Delta B=\bigl(A \setminus B\bigr)~\cup \bigl( B \setminus A \bigr)$. In a graph $G$, the \textit{distance} between two vertices $u$ and $v$, denoted $d(u,v)$, is the length of the shortest $uv$-path in $G$. The \textit{diameter} of $G$ is the length of a longest shortest path in $G$, i.e., $
\displaystyle{\underset{u,v \in V(G)}{\max} d(u,v)}$. 

Often, we will wish to modify one edge-coloring of some graph $G$ to obtain another edge-coloring. Formally, if $c$ is an edge-coloring of $G$ and $E' \subseteq E(G)$, we \textit{recolor} $E'$ by selecting edge-colors $c_i \neq c(e_i)$ for each $e_i \in E'$, and defining a new edge-coloring 
\[c'(e) = \begin{cases}  c(e) & \text{ if } e \not \in E', \\
c_i & \text{ if } e_i \in E'.
\end{cases}\]
In a slight abuse of notation, when we modify edge-colorings in this way, we will refer to both the original and the new edge-colorings as $c$. Given a recoloring of $E'$ and $e_i \in E'$, we say that $c(e_i)$ is a \textit{new color for c} if $c(e_i) \neq c(e)$ for all $e \not\in E'$.

Given a forbidden graph $F$, we will in particular wish to modify one rainbow $F$-free edge-coloring $c$ of $G$ to obtain another. We say $E' \subseteq E(G)$ is \textit{unrestricted} relative to a rainbow $F$-free edge-coloring $c$ if $G$ remains rainbow $F$-free under any recoloring of $E'$. Typically, $c$ and $F$ will be clear from context, and we will simply refer to such an edge set as unrestricted. In the case that $E' = \{e\}$, we simply say that the edge $e$ is unrestricted.

\section{Cliques}\label{sec:cliques}
\subsection{Bounds for general cliques}

To begin this section, we correct a proof from~\cite{Bushaw2022}, which states that for any $r \geq 4$, we have $\rsat(n, K_r) \geq (r - 1)n + O(1)$.

\begin{lem}\label{lem: KrLgDegSum}
    If $G$ is a rainbow $K_r$-saturated graph, then for any $u, v \in V(G)$ with $uv \not \in E(G)$, one of the following holds:
    \begin{enumerate}[(1)]
        \item $|N(u) \cap N(v)| \geq r - 1$, 
        \item $d(u) + d(v) \geq \binom{r}{2} - 1$.
    \end{enumerate}
\end{lem}
\begin{proof}
    Let $G$ be a rainbow $K_r$-saturated graph, $c$ be a proper edge-coloring of $G$ containing no rainbow copy of $K_r$, and $u, v \in V(G)$ with $uv \not \in E(G)$. Suppose $|N(u) \cap N(v)| < r - 1$. Note that the common neighborhood of any two nonadjacent vertices in a rainbow $K_r$-saturated graph has size at least $r - 2$, so we may assume $|N(u) \cap N(v)| = r - 2$. Let $G' = G[\{u, v\} \cup (N(u) \cap N(v))]$. Then, $G'$ must be a copy of $K_r^-$ (that is, $K_r$ with the edge $uv$ removed) which is rainbow under the coloring $c$, as otherwise, we may add the edge $uv$ to $G$ and extend $c$ without introducing a rainbow copy of $K_r$. Furthermore, if $|N(u) \triangle N(v)| < \binom{r - 2}{2}$, then there exists an edge $xy \in G[N(u) \cap N(v)]$ such that $c(xy)$ does not appear on any edges incident to $u$ or $v$. In this case, we may add the edge $uv$ to $G$ and assign it the color $c(xy)$, contradicting that $G$ is rainbow $K_r$-saturated. Therefore, it must be the case that $|N(u) \triangle N(v)| \geq \binom{r - 2}{2}$. Thus, we have
    \begin{equation*}
        \begin{split}
            d(u) + d(v)
                &=      |N(u) \triangle N(v)| + 2|N(u) \cap N(v)| \\
                &\geq   \binom{r - 2}{2} + 2(r - 2) \\
                &=      \binom{r}{2} - 1.
        \end{split}
    \end{equation*}
\end{proof}

\begin{lem}\label{lem:KrMinDeg}
    Let $r\geq 4$. If $G$ is a rainbow $K_r$-saturated graph, then $G$ has at most one vertex of degree $r - 2$.
\end{lem}

\begin{proof}
    Let $G$ be a rainbow $K_r$-saturated graph and $c$ be a proper edge-coloring of $G$ with no rainbow $K_r$. Suppose for contradiction that $G$ has two vertices $u, v$ of degree $r - 2$. 
    
    Note that if $uv \not \in E(G)$, then since $G$ is rainbow $K_r$-saturated, it must be the case that $N(u) = N(v)$ and $G[N(u) \cup \{u, v\}]$ is a rainbow copy of $K_{r}^-$ under $c$. But then, we may add the edge $uv$ to $G$ and color this edge with any color appearing in $G[N(u)]$ without creating a rainbow copy of $K_r$, a contradiction. 
    
    Therefore, we must have $uv \in E(G)$. However, in this case, adding any new edge to $u$ does not create a copy of $K_r$, a contradiction.
\end{proof}

\begin{prop}\label{prop: KrMinDegBd}
    Let $r\geq 4$ and $t \geq r - 1$. If $G$ is a rainbow $K_r$-saturated graph with $\delta(G) = t$, then
    \begin{equation*}
        e(G) \geq \left( \frac{r + t - 2}{2} \right)n + O(1).
    \end{equation*}
\end{prop}
\begin{proof}
    Let $G$ be a rainbow $K_r$-saturated graph with minimum degree $t \geq r - 1$. Let $u \in V(G)$ with $d(u) = t$. For any $v \in V(G) \setminus N[u]$, because $G$ is rainbow $K_r$-saturated, then $u$ and $v$ have at least $r - 2$ common neighbors and there exists a $K_{r - 2}$ in $G[N(u) \cap N(v)]$. Furthermore, $v$ has $d(v) - |N(u) \cap N(v)|$ neighbors outside of $N(u)$. Since this is true for any vertex not incident with $u$, we have:
    \begin{equation*}
        \begin{split}
            e(G)
                &=      e(N[u], V(G) \setminus N[u]) + e(G[V(G) \setminus N[u]]) + e(G[N[u]]) \\
                &\geq   (r - 2)(n - t) + \frac{1}{2}(t - (r - 2))(n - t - 1) + \binom{r - 2}{2} \\
                &=      \left( r - 2 + \frac{t - (r - 2)}{2} \right)n + O(1) \\
                &=      \left( \frac{t + r - 2}{2} \right)n + O(1).
        \end{split}
    \end{equation*}
\end{proof}

Together with Theorem~\ref{thm: k4 lower bound}, Proposition~\ref{prop: lower bound cliques} allows us to recover the rainbow saturation lower bound for cliques as stated in~\cite{Bushaw2022}.

\begin{prop}\label{prop: lower bound cliques}
    For $r \geq 5$, $sat^*(n, K_r) \geq (r - 1)n + O(1)$.
\end{prop}
\begin{proof}
    Let $G$ be a rainbow $K_r$-saturated graph. If $\delta(G) \geq r$, we get the desired bound by Proposition~\ref{prop: KrMinDegBd}. So, we may assume that $\delta:= \delta(G) \in \{r - 2, r - 1\}$. Let $u$ be a vertex of minimum degree and let $G' = G[V(G) \setminus N[u]]$. By Lemma~\ref{lem: KrLgDegSum}, for $v \not \in N[u]$, either $|N(u) \cap N(v)| \geq r - 1$ or $d(u) + d(v) \geq \binom{r}{2} - 1$.
    
    Let $T_1 = \{v \in V(G) \setminus N[u] : |N(u) \cap N(v)| \geq r - 1\}$ and let $T_2 = V(G) \setminus (N[u] \cup T_1)$. Then, for all $v \in T_2$, $d(v) \geq \binom{r}{2} - 1 - \delta$. Furthermore, since for all $v \in T_2$, $|N(u) \cap N(v)| = r - 2$, we obtain $d_{G'}(v) \geq \binom{r}{2} - 1 - \delta - (r - 2)$. Therefore, we have:
    \begin{equation*}
        \begin{split}
            e(G)
                &\geq   (r - 1)|T_1| + (r - 2)|T_2| + \frac{1}{2}|T_2|\left( \binom{r}{2} - 1 - \delta - (r - 2) \right) + O(1) \\
                &=      (r - 2)(|T_1| + |T_2|) + |T_1| + \frac{|T_2|}{2}\left( \binom{r}{2} - 1 - \delta - (r - 2) \right) + O(1) \\
                &\geq   (r - 2)(|T_1| + |T_2|) + |T_1| + \frac{|T_2|}{2}\left( \binom{r - 2}{2} - 1 \right) + O(1)\\
                &\geq   (r - 2)(|T_1| + |T_2|) + |T_1| + \frac{|T_2|}{2}\left( 2 \right) + O(1)\\
                &=      (r - 1)n + O(1).
        \end{split}
    \end{equation*}
\end{proof}

We give an upper bound construction. The key ingredient to this construction is Lemma~\ref{lem:rbwKrExt}. Given a graph $G$, a subgraph $H$ of $G$, an edge-coloring $c_H$ of $H$, and an edge-coloring $c_G$ of $G$. We say that $c_H$ \textit{extends} to $c_G$ if for any edge $e\in E(H)$ we have $c_{H}(e) = c_{G}(e)$, in which case we call $c_{G}$ the \textit{extension of $c_{H}$ to $G$}.

\begin{lem}\label{lem:rbwKrExt}
    Let $G = K_r$ and let $c: E(G) \to \mathbb{N}$ be a rainbow edge-coloring of $G$. Let $H = E_{r \binom{r - 1}{2} + 1}$ and set $G' := G \vee H$. Then, for any extension of $c$ to $G'$ that is a proper edge-coloring, there exists some vertex $v \in V(H)$ such that $G'[V(G) \cup \{v\}]$ is a rainbow $K_{r + 1}$ under $c$.
\end{lem}
\begin{proof}
    Let $G'$ be as described and let $c: E(G') \to \mathbb{N}$ be a proper edge-coloring of $G'$ such that $G'[V(G)]$ is a rainbow $K_r$. We will call a vertex $v \in V(H)$ bad if $G'[V(G) \cup \{v\}]$ is not a rainbow $K_{r + 1}$. For each $u \in V(G)$, there are at most $\binom{r - 1}{2}$ vertices in $V(H)$ that may repeat a color appearing on an edge of $G - u$. Since $|V(G)| = r$, then there are at most $r \binom{r - 1}{2}$ bad vertices of $H$. Since $|V(H)| \geq r \binom{r - 1}{2} + 1$, we are guaranteed to have at least one vertex which is not bad.
\end{proof}

The following construction allows us to obtain an explicit upper bound on the rainbow saturation for all cliques on at least four vertices. In a way, it refines the general upper bound proof in~\cite{Bushaw2022} when restricted to cliques. 

\begin{const}\label{kr sat construction}
    Let $r \geq 3$ and let $n \geq r \binom{r - 1}{2} + 2 + \sum_{i = 3}^{r} \left( i \binom{i - 1}{2} + 1 \right)$ and set $n' := n - 1 + \sum_{i = 3}^r \left( i \binom{i - 1}{2} + 1 \right)$. Now, let $G(r, n)$ be defined as follows:
    \begin{equation*}
        G(r, n) = E_1 \vee \left( \bigvee_{i = 3}^r E_{r\binom{r - 1}{2} + 1} \right) \vee E_{n'}.
    \end{equation*}

    We will call the $E_i$'s the parts of $G(r, n)$ and we will call the part of size $n'$ the leftover part of $G(r, n)$. Observe that $G(r, n)$ is a complete $r$-partite graph (and therefore, $K_{r + 1}$-free).
\end{const}
\begin{prop}
    Let $r \geq 3$ and $n \geq r \binom{r - 1}{2} + 2 + \sum_{i = 3}^{r} \left( i \binom{i - 1}{2} + 1 \right)$. Then, $G(r, n)$ as defined in Construction~\ref{kr sat construction} is properly rainbow $K_{r + 1}$-saturated.
\end{prop}
\begin{proof}
    We prove this by induction on $r$. To begin, let $r = 3$ and $n$ large enough. As mentioned before, since $G(3, n)$ is $3$-partite, then it is clearly rainbow $K_4$-free for any proper edge-coloring. Now, suppose we add an edge $e$ to $G(3, n)$ and let $c$ be any proper edge-coloring of $G(3, n) + e$. Because $G(3, n)$ is a complete $3$-partite graph, then $e$ must be contained within one of the parts of size greater than $1$. In this case, either $e$ is contained in the part of size $4$ or the leftover part, which has size at least $4$. In either case, $e$ is contained in a triangle with the part of size $1$. Finally, since any proper-edge coloring of a triangle leads to a rainbow triangle, then this triangle together with the part of size greater than $1$ not containing $e$ must contain a rainbow $K_4$ by Lemma~\ref{lem:rbwKrExt}. \\

    Next, suppose $r \geq 4$ and let $n$ be large enough (as defined in Construction~\ref{kr sat construction}). Suppose we have show that $G(r - 1, m)$ is properly rainbow $K_r$-saturated for all $m$ large enough. We claim that $G(r, n)$ is properly rainbow $K_{r + 1}$-saturated. As mentioned above, $G(r, n)$ is $r$-partite and therefore, clearly rainbow $K_r$-free for any proper edge-coloring. Denote the part of size $r \binom{r - 1}{2} + 1$ by $S_r$ and denote the leftover part by $S_{\ell}$. Let $n' = |S_{\ell}|$ and let $n_1 = n - |S_r|$, $n_2 = n - n'$. Observe that $G_1 := G(r,n)[V(G(r, n)) \setminus S_r] = G(r - 1, n_1)$ and $G_2 := G(r, n)[V(G(r, n)) \setminus S_{\ell}] = G(r - 1, n_2)$. Now, suppose we add an edge $uv$ to $G(r, n)$ and let $c$ be any proper edge coloring of $G(r, n) + uv$. Observe that the vertices $u$ and $v$ must be contained in at least one of $V(G_1)$ or $V(G_2)$. Without loss of generality, suppose $u$ and $v$ are contained in $V(G_1)$. Then, since $G_1 = G(r - 1, n_1)$, by induction, $G_1$ contains a rainbow copy of $K_r$, $K$. Now, by Lemma~\ref{lem:rbwKrExt}, $K \vee S_r$, a subgraph of $G(r, n)$ must contain a rainbow copy of $K_{r + 1}$. \\

    As $n \to \infty$, the number of edges of $G(r, n)$ not incident to a vertex in the leftover part is some constant. Therefore, we have:
    \begin{equation*}
        \begin{split}
            \lim_{n \to \infty} e(G(r, n))
                &= \left(\sum_{i = 2}^{r} r \binom{r}{2} + 1 \right)n + O(1) \\
                &= \left( \frac{r^4 - 2r^3 - r^2 + 10r - 8}{8} \right)n + O(1) \\
                &\leq \left( \frac{r^4}{8} \right) n + O(1).
        \end{split}
    \end{equation*}
\end{proof}

\begin{cor}
Let $r \geq 3$. Then, for $n$ large enough, we have:
\begin{equation*}
    \sat^{*}(n,K_{r + 1}) \leq \left( \frac{r^4}{8} \right)n + O(1).
\end{equation*}
\end{cor}

\subsection{Proof of Theorem~\ref{thm: k4 lower bound}}

In the case of $K_4$, we can obtain asymptotically tight bounds. Before we provide an upper bound construction for $\rsat(n,K_4)$, we make some observations which will be repeatedly useful.

\begin{obs}\label{shared nbhd edge}
Suppose $G$ is a properly rainbow $K_4$-saturated graph, and $xy \not \in E(G)$. Then $N(x) \cap N(y)$ contains an edge.

\begin{proof}
If $G$ is properly rainbow $K_4$-saturated, then $G$ has a rainbow $K_4$-free proper edge coloring $c$, but there is no rainbow $K_4$-free proper edge coloring of $G + xy$. In particular, we can properly edge-color $G + xy$ by coloring $E(G)$ according to $c$ and adding edge $xy$ in a color not appearing in $c$, so this proper edge-coloring of $G + xy$ contains a rainbow $K_4$-copy. Since $G$ is rainbow $K_4$-free under $c$, any rainbow $K_4$-copy in the described coloring must include $xy$. In particular, $xy$ is contained in some $K_4$-copy, say on $\{x,y,u,v\}$, and $uv$ is an edge in $N(x) \cap N(y)$.
\end{proof}

\end{obs}

\begin{obs}\label{K4 forbidden structure}

Suppose vertices $x,y,z$ form a $K_3$-copy, and let $V_j = \{v_1, v_2, \dots, v_j\}$ be a set of $j \geq 1$ vertices disjoint from $\{x,y,z\}$. Under any rainbow $K_4$-free proper edge-coloring of $\{x,y,z\} \vee V_j$, each $v_i \in V_j$ is adjacent to one of $x,y,z$ via an edge with color in $\{c(xy), c(xz), c(yz)\}$.
    
\end{obs}
\begin{proof}
Under any proper edge-coloring of $\{x,y,z\} \vee V_j$, the $K_3$ on $\{x,y,z\}$ must be rainbow. Moreover, for any $v_i \in V_j$, the edges $xv_i, yv_i, zv_i$ must receive distinct colors in any proper edge-coloring. Thus, if the $K_4$-copy on $\{x,y,z, v_i\}$ is not rainbow, then 
\[\{c(xy), c(xz), c(yz)\} \cap \{c(xv_i), c(yv_i), c(zv_i)\} \neq \emptyset.\]
\end{proof}

Observation~\ref{K4 forbidden structure} has two particular implications. Note that in a proper edge-coloring, $c(xy) \not\in \{c(xv_i), c(yv_i)\}$, so if $c(xy) \in \{c(xy), c(xz), c(yz)\} \cap \{c(xv_i), c(yv_i), c(zv_i)\}$, then $c(zv_i) = c(xy)$. Since $z$ is incident to at most one edge of color $c(xy)$, this implies that for at most one $v_i \in V_j$ do we have $c(xy) \in \{c(xy), c(xz), c(yz)\} \cap \{c(xv_i), c(yv_i), c(zv_i)\}$. Thus, as a consequence of Observation~\ref{K4 forbidden structure}, any properly rainbow $K_4$-saturated graph is $K_3 \vee E_4$ free. Moreover, if $G$ is a properly edge-colored, rainbow $K_4$-free graph containing a copy of $K_3 \vee E_3$, then the edge-colors used in this $K_3$-copy must be repeated on a matching between $V(K_3)$ and $V(E_3)$. Using Observation~\ref{K4 forbidden structure}, we can now quickly show that Construction~\ref{k4 sat construction} is properly rainbow $K_4$-saturated. This construction will provide the upper bound in Theorem~\ref{thm: k4 lower bound}.

\begin{const}\label{k4 sat construction}

For $n\geq 6$, consider the $n$-vertex graph $G(n)$ obtained as follows. Let $G'(n)$ be the $n-2$ vertex graph on $k := \lceil \frac{n-2}{4} \rceil$ components $C_1, C_2, \dots, C_k$, with $C_i = K_4$ if $i \leq \lfloor \frac{n-2}{4} \rfloor$. We label the vertices of $C_i$ as $\{v_{i,j}: 1 \leq j \leq 4\}$. If $n-2 \equiv m \mod 4$ with $m \neq 0$, then let $C_k = K_{m}$, and label the vertices of $C_k$ as $\{v_{k,j}: 1\leq j \leq m\}$. We take $G(n) = K_2 \vee G'(n)$, labeling the vertices of the $K_2$-copy joined to $G'(n)$ as $x,y$.

We edge-color $G(n)$ as follows. Color $G(n)[C_1 \cup \{x,y\}]$ with five colors using the perfect matching decomposition of a $K_6$. Suppose $c(xy) = 0$ in this coloring and the other colors used are $\{1, 2, 3, 4\}$. For $1 < i \leq \left \lfloor \frac{k - 2}{4} \right \rfloor$, color the $K_6$ induced by $C_i \cup \{x, y\}$ the exact same (most importantly, ensuring $c(xy) = 0$), switching colors $\{1, 2, 3, 4\}$ with colors $\{4i + 1, 4i + 2, 4i + 3, 4i + 4\}$. Finally, if $|V(C_k)| < 4$, we color the edges incident to $V(C_k)$ in any legal fashion (using at most $4$ new colors) which avoids a rainbow $K_4$-copy on $V(C_k) \cup \{x,y\}$.
\end{const}

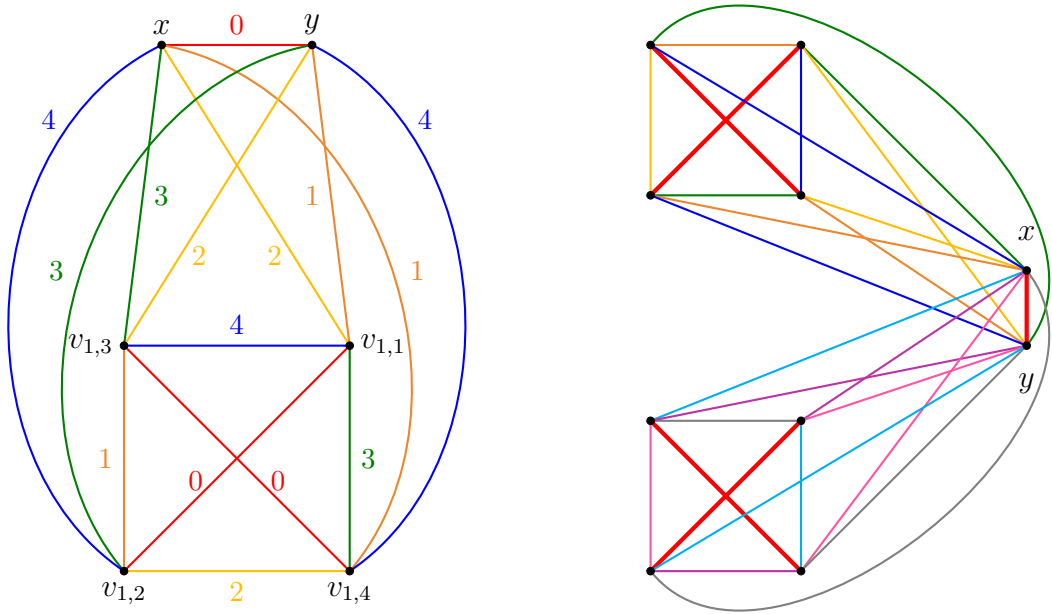
\begin{figure}[h]

\begin{center}

\begin{tikzpicture}

\draw[thick,red] (-4.5,7) -- (-6.5,7) node[pos = 0.5, above]{\small 0};
\draw[thick,red] (-4,0) -- (-7,3) node[pos = 0.4, right]{\small 0};
\draw[thick,red] (-4,3) -- (-7,0) node[pos = 0.6, left]{\small 0};

\draw[thick, cadmium] (-7,3) -- (-7,0) node[pos = 0.5, left]{\small 1};
\draw[thick, cadmium] (-4,3) -- (-4.5,7) node[pos = 0.5, left]{\small 1};
\draw[thick, cadmium] (-4,0) to[bend right = 60] (-6.5,7);
\draw[cadmium] (-3.1,4) node{\small 1};

\draw[thick, amber] (-7,0) -- (-4,0) node[pos = 0.5, below]{\small 2};
\draw[thick, amber] (-4,3) -- (-6.5,7) node[pos = 0.3, left]{\small 2};
\draw[thick, amber] (-7,3) -- (-4.5,7) node[pos = 0.3, right]{\small 2};

\draw[thick, forest] (-4,0) -- (-4,3) node[pos = 0.5, right]{\small 3};
\draw[thick, forest] (-7,0) to[bend left = 60] (-4.5,7);
\draw[forest] (-7.9,4) node{\small 3};
\draw[thick, forest] (-7,3) -- (-6.5,7) node[pos = 0.5, right]{\small 3};

\draw[thick, blue] (-4,3) -- (-7,3) node[pos = 0.5, above]{\small 4};
\draw[thick, blue] (-4,0) to[bend right = 60] (-4.5,7);
\draw[blue] (-8,6) node{\small 4};
\draw[thick, blue] (-7,0) to[bend left = 60] (-6.5,7);
\draw[blue] (-3,6) node{\small 4};

\filldraw (-4,0) circle (0.05 cm) node[below]{$v_{1,4}$};
\filldraw (-4,3) circle (0.05 cm) node[right]{$v_{1,1}$};
\filldraw (-7,0) circle (0.05 cm) node[below]{$v_{1,2}$};
\filldraw (-7,3) circle (0.05 cm) node[left]{$v_{1,3}$};

\filldraw (-4.5,7) circle (0.05 cm) node[above]{$y$};
\filldraw (-6.5,7) circle (0.05 cm) node[above]{$x$};

\draw[thick, cadmium] (0,7) -- (2,7); 
\draw[ultra thick, red] (0,7) -- (2,5); 
\draw[thick, amber] (0,7) -- (0,5); 
\draw[thick, forest] (0,7) to[bend left = 90] (5,3); 

\draw[thick, blue] (2,7) -- (2,5); 
\draw[ultra thick, red] (2,7) -- (0,5); 
\draw[thick, forest] (2,7) -- (5,4); 
\draw[thick, amber] (2,7) -- (5,3); 

\draw[thick, forest] (2,5) -- (0,5); 
\draw[thick, amber] (2,5) -- (5,4); 
\draw[thick, cadmium] (2,5) -- (5,3); 

\draw[thick, cadmium] (0,5) -- (5,4); 
\draw[thick, blue] (0,5) -- (5,3); 

\draw[thick, blue] (0,7) -- (5,4); 

\draw[thick, gray] (0,2) -- (2,2); 
\draw[ultra thick, red] (0,2) -- (2,0); 
\draw[thick, brilliantrose] (0,2) -- (0,0); 
\draw[thick, cyan] (0,2) -- (5,4); 
\draw[thick, byzantine] (0,2) -- (5,3); 

\draw[thick, cyan] (2,2) -- (2,0); 
\draw[ultra thick, red] (2,2) -- (0,0); 
\draw[thick, byzantine] (2,2) -- (5,4); 
\draw[thick, brilliantrose] (2,2) -- (5,3); 

\draw[thick, byzantine] (2,0) -- (0,0); 
\draw[thick, brilliantrose] (2,0) -- (5,4); 
\draw[thick, gray] (2,0) -- (5,3); 

\draw[thick, gray] (0,0) to[bend right = 90] (5,4); 
\draw[thick, cyan] (0,0) -- (5,3); 

\draw[ultra thick, red] (5,4) -- (5,3); 

\draw[fill = black] (0,2) circle (0.05 cm);
\draw[fill = black] (2,2) circle (0.05 cm);
\draw[fill = black] (2,0) circle (0.05 cm);
\draw[fill = black] (0,0) circle (0.05 cm);

\draw[fill = black] (0,7) circle (0.05 cm);
\draw[fill = black] (2,7) circle (0.05 cm);
\draw[fill = black] (2,5) circle (0.05 cm);
\draw[fill = black] (0,5) circle (0.05 cm);

\draw[fill = black] (5,4) circle (0.05 cm);
\draw[fill = black] (5,3) circle (0.05 cm);

\node at (5,4.5) {$x$};
\node at (5,2.5) {$y$};

\end{tikzpicture}

\caption{Construction~\ref{k4 sat construction} on 6 and 10 vertices; note the repetition of color 0 when $x,y$ are adjacent to multiple $K_4$-copies.}\label{k4 illustration}

\end{center}

\end{figure}

We illustrate the described coloring of an edge $xy$ and the components $C_1$ and $C_2$ in Figure~\ref{k4 illustration}. The proof of Proposition~\ref{upper bound sat} yields the upper bound of Theorem~\ref{thm: k4 lower bound}.

\begin{prop}\label{upper bound sat}
Construction~\ref{k4 sat construction} is properly rainbow $K_4$-saturated for all $n \geq 6$.
\end{prop}
\begin{proof}
    Let $c$ be the edge-coloring of $G(n)$ described in Construction~\ref{k4 sat construction}. Since there are no edges between any $C_i$, then any possible rainbow $K_4$ under $c$ must be contained in the graph induced by $C_i \cup \{x, y\}$ for some $1< i\leq \lfloor \frac{k-2}{4} \rfloor$. However, for any $i$, $G(n)[C_i \cup \{x, y\}]$ contains at most $5$ colors, ensuring that no rainbow $K_4$ can exist. \\
    Now, suppose we add an edge $v_{i, j}v_{i',j'}$ to $G(n)$. Without loss of generality, suppose $i \leq \left \lfloor \frac{n - 2}{4} \right \rfloor$. In particular, $|C_i| = 4$. Then, observe that $v_{i,j}, x, y$ form a triangle in $G(n)$ and the vertices $v_{i',j'}, v_{i, j + 1}, v_{i, j + 2}, v_{i, j + 3}$ (where the second indices are taken mod $4$) are joined to every vertex of the triangle by an edge in $G(n) + v_{i, j}v_{i',j'}$. That is, $G(n) + v_{i, j}v_{i',j'}$ contains a copy of $K_3 \vee E_4$. Therefore, by Observation~\ref{K4 forbidden structure}, $G(n) + v_{i, j}v_{i',j'}$ contains a rainbow $K_4$ under any coloring.
\end{proof}

We now turn to proving the lower bound in Theorem~\ref{thm: k4 lower bound}. We begin with a lemma which greatly restricts the number and behavior of very low-degree vertices in a properly rainbow $K_4$-saturated graph.

\begin{lem}\label{degree 3 bound}
If $G$ is rainbow $K_4$-saturated, then the vertices of $G$ with degree at most $3$ form a clique.
\end{lem}

\begin{proof}
Let $u,v$ be vertices with degree at most $3$, and suppose $uv \not \in E(G)$. By Observation~\ref{shared nbhd edge}, there is an edge $xy$ in $N(u)\cap N(v)$. 
Without loss of generality, $2\leq d(u)\leq d(v) \leq 3$. Fix a rainbow $K_4$-free coloring of $G$, say $c$. 
There are four cases to consider, ordered by complexity. 

\begin{enumerate}
        \item[Case 1:] $d(u) = d(v) = 2$.
        
        Observe that the edges $ux$ and $vy$ are not contained in any $K_4$ copy in $G$, so $\{ux, vy\}$ is unrestricted. We may thus re-color so that $c(ux) = c(vy)$ is a new color for $c$. Now, adding the edge $uv$ in another new color creates a proper edge-coloring of $G + uv$. Note that $uv$ is contained in exactly one $K_4$ copy, on $\{u,v,x,y\}$, which is not rainbow. Thus, $G + uv$ admits a rainbow $K_4$-free proper edge-coloring, a contradiction.
                
        \item[Case 2:] $d(u) = 2, d(v)=3$.
        
        Let $z\in N(v)\setminus N(u)$. Since $uz\notin E(G)$, there must exist an edge in $N(u) \cap N(z)$, which must be $xy$ since $N(u) = \{x,y\}$. Hence $\{x,y,z\}$ form a clique. The edges $ux,uy$ are in no $K_4$, so $\{ux,uy\}$ is unrestricted. Both $vx$ and $vy$ are in exactly one $K_4$, on vertex set $\{v,x,y,z\}$. It follows that either $\{ux, uy, vx\}$ or $\{ux, uy, vy\}$ is unrestricted: if $c(xy)=c(vz)$, then both are unrestricted; if $c(xz) = c(vy)$, then $\{ux, uy, vx\}$ is unrestricted; if $c(yz) = c(vx)$, then $\{ux, uy, vy\}$ is unrestricted. Without loss of generality, $\{ux, uy, vy\}$ is unrestricted. By recoloring $\{ux,vy\}$ as in Case $1$, we can extend $c$ to a rainbow $K_4$-free edge-coloring of $G + uv$, a contradiction.

        \item[Case 3: ]  $d(u) = d(v) = 3$, $N(u)\neq N(v)$.
        
        Let $z\in N(u)\setminus N(v)$ and $w\in N(v)\setminus N(u)$. First, suppose $c(uz)\neq c(xy)$ and $c(vw)\neq c(xy)$. We extend $c$ to $G + uv$ by setting  $c(uv) = c(xy)$. This extension of $c$ remains proper, and $uv$ is contained in one $K_4$ copy, on $\{u,v,x,y\}$, which is not rainbow, a contradiction. Hence either $c(uz)=c(xy)$ or $c(vw) = c(xy)$. Relabeling if necessary, we may assume $c(uz) = c(xy) = 1$. Then $\{ux, uy\}$ is unrestricted. If $\{x,y,w\}$ does not induce a clique, then in fact $\{ux, uy, vx, vy \}$ is unrestricted. If so, we may proceed as in Case $1$ to reach a contradiction. Hence we may assume $\{x,y,w\}$ induces a clique. As in Case $2$, at least one of $\{ux, uy, vx\}$ and $\{ux, uy, vy\}$ is free. Without loss of generality, $\{ux, uy, vy\}$ is free, and we recolor $ux$ and $vy$ with the same new color, as in previous cases. Adding the edge $uv$ in any color introduces no rainbow $K_4$, a contradiction.
        
        \item[Case 4:] $d(u) = d(v) = 3$, $N(u) = N(v)$.
        
        Say $N(u) = \{x,y,z\}$. If $\{x,y,z\}$ does \textit{not} induce a clique, then the set of all edges incident to $u$ and $v$ is unrestricted. Then we may recolor so that $c(ux) = c(vy) = k_1$, $c(uy) = c(vz) = k_2$, $c(uz) = c(vx)=k_3$, where all $k_i$ are new colors. We then extend $c$ to $G + uv$ letting $c(uv)$ be another new color.  The only copies of $K_4$ that can contain $uv$ in $G + uv$ are on vertex sets $\{u,x,y,v\}$, $\{u,x,z,v\}$, and $\{u,y,z,v\}$, none of which are rainbow under $c$. Hence $G+uv$ is rainbow $K_4$-free under $c$, a contradiction. \\
        Hence, $\{x,y,z\}$ induces a clique. Then $\{u,x,y,z\}$ is a $K_4$ copy, and so two of its edges must receive the same color under $c$. Without loss of generality, $c(xy) = c(uz) = 1$. Similarly, two edges induced by $\{v,x,y,z\}$ must repeat a color under $c$. This can not be the same color repeated in $\{u,x,y,z\}$. Without loss of generality, $c(yz) = c(vx) = 2$. Note $c(xz)$ must be distinct from both of these, say $c(xz) = 3$. Then $\{ux,uy,vy,vz\}$ is unrestricted. We recolor so that $c(ux) = c(vy) = k_1$, $c(uy) = c(vz)=k_2$, where $k_i$ are new colors. We add the edge $uv$, and extend $c$ to $G + uv$ by setting $c(uv) = 3$. Note, $c$ remains proper since $u$ and $v$ are incident to no edges colored $3$ under $c$ in $G$. Adding the edge $uv$ creates three new $K_4$ copies: the copy on $\{u,x,y,v\}$ repeats color $k_1$, the copy on $\{u,x,z,v\}$ repeats color $3$, and the copy on  $\{u,y,z,v\}$ repeats color $k_2$. Hence $G+uv$ admits a rainbow $K_4$-free proper edge-coloring, a contradiction.
\end{enumerate}
    In any case, we reach a contradiction. Hence, $uv\in E(G)$. This holds for any pair of vertices of degree at most $3$. Hence all vertices of degree at most $3$ form a clique.
\end{proof}

We now derive the main tool (Lemma~\ref{nice dominating set}) needed to prove Theorem~\ref{thm: k4 lower bound}. The broad strategy of our argument is to first argue that any rainbow $K_4$-saturated graph $G$ has a dominating set $D$ with certain useful properties, and then count edges in $G$ by examining the structure of $G \setminus D$ and the density of edges between $D$ and $V(G \setminus D)$. Towards this end, we begin by arguing that any sparse rainbow $K_4$-saturated graph must contain a small dominating set satisfying one of two properties. We say that $D \subset V(G)$ is \textit{k-dominating} if for any $v \notin D$, we have $|N(v) \cap D| \geq k$.

\begin{lem}\label{nice dominating set}
Fix $\frac{1}{2} > \alpha > 0$, and suppose $G$ is an $n$-vertex, rainbow $K_4$-saturated graph with $n$ large enough that $\alpha^2 n \geq 7$. Then either $e(G) \geq \frac{7}{2}n - 4 \alpha n$, or $G$ contains a dominating set $D$ with $|D| \leq 2 \alpha n$ such that $D$ satisfies one of the following conditions:

\begin{enumerate}

\item $D$ is a $3$-dominating set;

\item $D$ is a $2$-dominating set containing adjacent vertices $x,y$ such that $v \in N(x,y)$ for every $v \in V(G)\setminus D$.

\end{enumerate}

\end{lem}

\begin{proof}

Let $v$ be a vertex of minimum degree in $G$. If $d(v) \geq \alpha n - 1$, then by choice of $n$, we have $\delta(G) \geq \frac{\alpha}{2}n$, and thus
\[e(G) \geq \frac{1}{2} \cdot \frac{\alpha}{2} n^2 = \frac{\alpha}{4}n^2 > \frac{\alpha^2}{2}n^2,\]
with the last inequality following from $\frac{1}{2} > \alpha$. Now, since $\alpha^2 n \geq 7$, we have 
\[e(G) \geq \frac{\alpha^2}{2}n^2 \geq \frac{7}{2}n,\] 
and so the desired edge bound holds.

Thus, we may assume that $|N[v]| \leq \alpha n$. If $N[v]$ is a $3$-dominating set, then we set $D = N[v]$ and are done. Suppose not. We let $V' := V(G) \setminus N[v]$, and $G' := G[V']$. By Observation~\ref{shared nbhd edge}, $|N(u) \cap N(v)| \geq 2$ for any $u \in V'$, so $N[v]$ is a $2$-dominating set. We now consider the degrees of vertices in $V'$, and in particular those vertices which have either two neighbors in $N[v]$ or few neighbors in $V'$.

We define
\[L_1 := \{u \in V': |N(u) \cap N(v)| = 2\}\]
and
\[L_2 := \{u \in V' : |N(u) \setminus N(v)| \leq 2 \}. \]

We also define a set of vertices with many neighbors in $N[v]$:

\[H := \{u \in V': |N(u) \cap N(v)| \geq 4 \}.\]
Note that $H$ and $L_1$ are disjoint.

We first show that we may assume $L_1 \cap L_2$ is not empty. Indeed, if $L_1 \cap L_2 = \emptyset$, then we can bound

\begin{equation}\label{cross edge bound} 
e(V', N[v]) \geq 4|H| +  3(|V'| - |H| ) - |L_1| = 3|V' \setminus (L_1 \cup L_2)| + 2|L_1| + 3|L_2| + |H|
\end{equation}
and estimate $e(G')$ as follows. Recall that by Lemma~\ref{degree 3 bound}, the vertices of degree at most $3$ in $G$ form a clique. In particular, Lemma~\ref{degree 3 bound} implies that $L_2 \setminus H$ contains at most one vertex $u$ with $d_{G'}(u) = 0$. Indeed, if $u_1, u_2 \in L_2 \setminus H$ have degree $0$ in $G'$, then $d(u_1), d(u_2)$ are both at most $3$, which implies $u_1u_2$ is an edge in $G'$. Thus, at most one vertex of $L_2 \setminus H$ has degree $0$ in $G'$, so we have
\begin{align}\label{G' edge bound}
e(G') = \frac{\sum_{w \in V' \setminus L_2} d(w) +
\sum_{u \in L_2} d_{G'}(u)}{2} \geq \frac{3|V' \setminus L_2| + |L_2 \setminus H| - 1}{2} \nonumber \\
= \frac{3}{2}|V' \setminus (L_1 \cup L_2)| + \frac{3}{2}|L_1| + \frac{1}{2}|L_2 \setminus H| - \frac{1}{2}.
\end{align}
We combine (\ref{cross edge bound}) and (\ref{G' edge bound}), noting that since $(L_2\setminus H) \cup H \supseteq L_2$, we have $3|L_2| + |H| + \frac{1}{2}|L_2 \setminus H| \geq \frac{7}{2}|L_2|$. Thus,
\[e(G) \geq \frac{9}{2}|V \setminus (L_1 \cup L_2)| + \frac{7}{2}|L_1| + \frac{7}{2}|L_2| - \frac{1}{2} \geq \frac{7}{2}|V'| - 1 \geq \frac{7}{2}n - \frac{7\alpha}{2}n - 1 > \frac{7}{2}n - 4\alpha n,\]
yielding the desired edge bound.

So we may assume that $L_1 \cap L_2$ is not empty; let $z \in L_1 \cap L_2$, with $N(v) \cap N(z) = \{x,y\}$. Since $N(v) \cap N(z)$ is not an independent set, we must have $xy \in E(G)$. Our goal now is to show that almost all vertices in $V'$ are in $N(x,y)$. If $d(z) = 2$, then by Observation~\ref{shared nbhd edge}, any $u \in V'$ is in $N(x,y)$. So, assume $d(z) \geq 3$. Suppose there exists $w \in V' \setminus [N(x,y) \cup N(z)]$. (So, $w$ is not dominated by $xy$, and $w$ is not adjacent to $z$.)

\begin{claim}\label{augmented dominating set}

If $|N(w)| \leq \alpha n$, then $N[v] \cup N[w] \cup N(z)$ is a $3$-dominating set of size at most $2 \alpha n$.

\end{claim}

\begin{proof}[Proof of Claim~\ref{augmented dominating set}]

We set $D := N[v] \cup N[w] \cup N(z)$, and begin by bounding $|D|$. Since $z \in L_1 \cap L_2$, we have $|N(z) \cap V'| \leq 2$. Since $zw \notin E(G)$, we must have $|N(z) \cap N(w)| \geq 2$; in particular, since $N(w)$ contains at most one of $x,y$, we know that $|N(z) \cap V'|$ is non-empty and intersects $N(w)$. Thus, $|N[w] \cup (N(z) \cap V')| \leq |N[w]| + 1 \leq \alpha n + 2$. Moreover, since $vw \notin E(G)$, we must have $|N[w] \cap N[v]| \geq 2$, implying that
\[ |D| = |N[v]| + |(N[w] \cup N(z)) \cap V'| \leq |N[v]| + |N[w] \cup (N(z) \cap V')| - 2 \leq 2\alpha n,\]
as desired. 

Next, we show that every vertex in $V \setminus D$ has at least $3$ neighbors in $D$. By construction, $z$ has at least $3$ neighbors in $D$ since $|N(z)| \geq 3$. Take $u \neq z$ in $V \setminus D$. Thus, $u$ is not adjacent to $v, w,$ or $z$. In particular, $u$ must have at least two neighbors in each of $N(v), N(w),$ and $N(z)$. We have two cases.

\begin{enumerate}

\item $x,y \in N(u)$.

In this case, since $|\{x,y\} \cap N(w)| \leq 1$, we know that $u$ must have at least one more neighbor in $N(w)$. Thus, $|N(u) \cap D| \geq 3$.

\item $N(u)$ contains at most one of $x,y$.

In this case, we know that $|N(u) \cap N(v)| 
\geq 2$ and that $N(u)$ contains a neighbor of $z$ which is not in $\{x,y\}$. Since $N(v) \cap N(z) = \{x,y\}$, this implies that $N(u)$ intersects $D \setminus N(v)$;  thus, $|N(u) \cap D| 
\geq 3$.
\end{enumerate}

Hence, $D$ is a $3$-dominating set, completing the claim.
\end{proof}

Now, by Claim~\ref{augmented dominating set}, we are done if there exists a vertex $w \in V'\setminus [N(x,y) \cup N(z)]$ with $d(w) \leq \alpha n$. On the other hand, if $V$ contains at least $\alpha n - 8$ vertices of degree greater than $\alpha n$, then 
\[e(G) \geq \frac{1}{2} (\alpha n)(\alpha n - 8) = \frac{\alpha^2}{2}n^2 - 4\alpha n \geq \frac{7}{2}n - 4\alpha n,\]
satisfying the desired bound. So $|V' \setminus [N(x,y) \cup N(z)]| < \alpha n - 8$, and thus 
\[D := N[v] \cup [V' \setminus N(x,y)]\]
has $|D| \leq 2 \alpha n$, and $D$ is a $2$-dominating set with the property that $u \in N(x,y)$ for any $u \notin D$. 
\end{proof}

We are now ready to prove Theorem~\ref{thm: k4 lower bound}, which we restate here for convenience.

\kfourlowerbound*

\begin{proof}

By Proposition~\ref{upper bound sat}, we have $\rsat(n,K_4) \leq \frac{7}{2}n + O(1)$. Thus, we must only demonstrate the lower bound. Suppose that $G$ is an $n$-vertex, properly rainbow $K_4$-saturated graph with $\alpha^2 n \geq 7$. If $e(G) \geq \frac{7}{2}n - 4 \alpha n$, we are done, so suppose not. By Lemma~\ref{nice dominating set}, $G$ contains a dominating set $D$ with $|D| \leq 2 \alpha n$ satisfying one of the outcomes of Lemma~\ref{nice dominating set}.
We now have two cases, depending upon which outcome occurs. For notational convenience, we set $V':=V(G) \setminus D$ and let $G'$ be the subgraph of $G$ induced on $V'$. 

Suppose first that $D$ is a $3$-dominating set. Recall that by Lemma~\ref{degree 3 bound}, the set of vertices in $G$ of degree at most $3$ form a clique. Note that if $v \in V'$ has degree at most $3$, then in fact $d(v) = 3$ and $N(v) \subseteq D$. Thus, the set of vertices of degree at most $3$ intersects $V'$ at most once. We define 

\[V_1 := \{v \in V' : d_{G'}(v) = 0\}\]
and 
\[V_2 := \{v \in V' : d_{G'}(v) \geq 1\}.\]
By the above observations, $V_1$ contains at most one vertex of degree $3$ in $G$. We have

\begin{align*}
e(G) &\geq e(V', D) + e(G') = \sum_{v \in V_1} d(v) + \sum_{u \in V_2} |N(u) \cap D| + \frac{1}{2} \sum_{u \in V_2} d_{G'}(u)& \\
&\geq 4|V_1| - 1 + 3|V_2| + \frac{|V_2|}{2} \geq \frac{7}{2} |V'| - 1.&
\end{align*}

Since $|V'| = n - |D|  \geq (1 - 2\alpha)n$, we have 
\[e(G) \geq \frac{7}{2}(1 - 2 \alpha)n - 1 = \frac{7}{2}n - 7 \alpha n - 1 > \frac{7}{2}n - 8 \alpha n,\]
as desired.

Next, suppose $D$ is a $2$-dominating set containing adjacent vertices $x,y$ such that $v \in N(x,y)$ for every $v \in V(G)\setminus D$. In this case, note that $V(G) \setminus N(x,y) \subseteq D$, and both are $2$-dominating sets of the desired type. For later convenience, we will redefine $D$ if necessary, setting
\[D := V(G) \setminus N(x,y).\]
Note that now, $V' = V(G)\setminus D$ is simply equal to $N(x,y)$. As before, we set $G' = G[V']$.

Now, we begin by making several observations on the structure and coloring of $G'$. In particular, we will find a matching within $G'$ of edges which receive the same color as $xy$; to indicate the special role played by this edge-color, we shall say $c(xy) = 0$. We begin by showing that certain vertices in $G'$ must be incident to edges of color $0$.

Suppose $z \in V'$ has $d_{G'}(z) \geq 3$. By Observation~\ref{K4 forbidden structure}, we must have $d_{G'}(z) = 3$, and each of the three vertices in $N(z) \cap V'$ is adjacent to one of $x,y,z$ via an edge with color in $\{0, c(xz), c(yz)\}$. To maintain a proper edge-coloring, only $x$ can be adjacent to $N(z) \cap V'$ via an edge of color $c(yz)$, only $y$ can be adjacent to $N(z) \cap V'$ via an edge of color $c(xz)$, and only $z$ can be adjacent to $N(z) \cap V'$ via an edge of color $0$. Thus, the maximum degree in $G'$ is $3$, and every vertex of degree $3$ in $G'$ is incident to an edge of color $0$.

We now investigate vertices with degree $2$ in $G'$. In particular, we define the set of vertices with degree $2$ in $G'$ and exactly $2$ neighbors in $D$, 
\[S_2^2 := \{v \in V': d_{G'}(v) = 2 \text{ and } |N(v)\cap D| = 2\}.\]
Note that every $v \in S_2^2$ has $N(v) \cap D = \{x,y\}$. We start by noting that $S_2^2$ contains very few vertices which are not incident to an edge of color $0$. 

Indeed, suppose there exist $u,v \in S_2^2$ such that $d(u,v) > 2$ and neither $u$ nor $v$ is incident to an edge of color $0$. We can thus add $uv$ to $G$ with $c(uv) = 0$ while maintaining a proper coloring. Moreover, since $d(u,v) > 2$, the addition of $uv$ does not form a triangle in $G'$. Thus, any copy of $K_4$ using $uv$ in $G + uv$ must contain two vertices from $D$. Since $u,v \in S_2^2$, their only neighbors in $D$ are $x,y$. Hence, $\{u,v,x,y\} $ span the unique copy of $K_4$ using $uv$ in $G + uv$; since $c(uv) = c(xy) = 0$, this copy is not rainbow. So $G + uv$ can be properly edge-colored while avoiding a rainbow $K_4$-copy, a contradiction to the assumption that $G$ is rainbow $K_4$-saturated. 

Thus, the subset of vertices in $S_2^2$ which are not incident to edges of color $0$ pairwise are at distance at most $2$ in $G'$. Since $\Delta(G') \leq 3$, for \textit{any} vertex $u \in G'$, we have
\[ |\{ v \in G' : d_{G'}(u,v) \leq 2 \}| \leq 13. \]
While a tighter bound in fact holds for vertices in $S_2^2$, we shall only require that the subset of vertices in $S_2^2$ not incident to edges of color $0$ is of constant size.

Now, let $M_0$ be the matching in $G'$ consisting of all edges in $E(G')$ colored 0. We form $G'_0$ by deleting $M_0$ from $G'$. Now, by the above observations, $\Delta(G'_0) \leq 2$. Thus, every component of $G'_0$ is a path (possibly trivial) or a cycle. We shall use this structure to count edges in $G$. 

Call $v \in V'$ \textit{light} if $N(v) \cap D = \{x,y\}$, and \textit{heavy} if not. (Thus, all $S_2^2$ vertices are light, but $V'$ may contain light vertices which are not in $S_2^2$.) Call a component $C$ of $G'_0$ \textit{light} if at most one vertex of $C$ is heavy, and let $L$ be the set of light components of $G'_0$. 

Similarly, call a component $C$ \textit{heavy} if at least two vertices of $C$ are heavy, with the set of heavy components denoted by $H$. For a component $C$ of $G'_0$, let $L_C$ denote the set of light vertices in $V(C)$ and $H_C$ the set of heavy vertices in $V(C)$.  

Now, $e(G) \geq e(G',D) + e(G')$, which we can bound by summing component-by-component:  

\begin{align*} 
e(G) \geq &\sum_{v \in V'} \left(|N(v) \cap D| + \frac{1}{2}d_{G}'(v)\right) = \sum_{v \in C, C \in L} \left(|N(v) \cap D| + \frac{1}{2}d_{G'}(v)\right)& \\ & + \sum_{v \in C, C \in H} \left(|N(v) \cap D| + \frac{1}{2}d_{G'}(v)\right).& 
\end{align*}

In particular, if the number of edges incident to component $C$ satisfies 
\begin{equation}\label{good component density}
i(C) =: \sum_{v \in C} \left(|N(v) \cap D| + \frac{1}{2}d_{G'}(v)\right) \geq \frac{7}{2}|C|
\end{equation}
for all components $C$, then the desired bound on $e(G)$ would hold. We will not be able to guarantee that (\ref{good component density}) holds for all $C$, but show that only constantly many components $C$ can fail to satisfy (\ref{good component density}). 

We shall obtain different bounds on $i(C)$ for different types of components. Firstly, recall that by Lemma~\ref{degree 3 bound}, at most four vertices of $G$ have degree $\leq 3$, and we have argued above that at most $13$ vertices in $S_2^2$ are not incident to an edge of color $0$ in $G'$. Call a component $C$ \textit{exceptional} if $C$ contains either a vertex with degree $\leq 3$ in $G$ or a vertex in $S_2^2$ which is not incident to an edge of color $0$ in $G'$. Note that there are at most $17$ exceptional components of $G'_0$ in total. We now estimate $i(C)$ for different types of components $C$.

We first address components of size $1$.   
If $C$ is a component such that some vertex $v$ in $C$ has $d_{C}(v) = 0$, then $V(C) = \{v\}$. In this case, $i(C) = d(v)$. Note that if $C$ is non-exceptional, then $d(v) \geq 4$, so $i(C) > \frac{7}{2}|C|$.

Next, we consider components of size greater than $1$. We treat these in cases, depending upon their classification as light or heavy.

\begin{enumerate}

\item[Case 1:] $|V(C)|\geq 2$ and $C \in H$.

For a non-exceptional component $C \in H$, we estimate $i(C)$ as follows. 
If $v$ is a light vertex of $C$, then (since $C$ is non-exceptional) $v$ either is in $S_2^2$ and is incident to a color 0 edge in $G'$ or $d_{G'}(v) = 3$. In either case, $v$ is incident to an edge in $M_0$, so $d_{G'}(v) = d_{C}(v) + 1$. 
If $v$ is a heavy vertex of $C$, then $|N(v) \cap D| \geq 3$, and (since $|V(C)| \geq 2$) $d_{C}(v) \geq 1$. Now,
\begin{align*}i(C) &= \sum_{u \in L_C} \left( |N(u) \cap D| + \frac{d_{G'}(u)}{2}\right) + \sum_{v \in H_C} \left( |N(v) \cap D| + \frac{d_{G'}(v)}{2}\right)& \\ &\geq \sum_{u \in L_C} \left(2 + \frac{d_C(u) + 1}{2}\right) + \sum_{c \in H_C} \left(3 + \frac{d_{C}(v)}{2} \right).& 
\end{align*}
Since $C$ is either a path or a cycle, at most two vertices of $C$ have degree $1$ in $C$, so 
\[\sum_{u \in L_C} \left(2 + \frac{d_C(u) + 1}{2}\right) + \sum_{c \in H_C} \left(3 + \frac{d_{C}(v)}{2} \right)\] \[ \geq \frac{7}{2}|L_C| + 4|H_C| - 1 = \frac{7}{2}|C| + \frac{1}{2}|H_C| - 1 \geq \frac{7}{2}|C|\]
since, because $C$ is heavy, $|H_C| \geq 2$.

If $C$ is an exceptional heavy component, then at most 17 vertices in $C$ are either of degree at most $3$ in $G$ or are in $S_2^2$ but not incident to an edge of $M_0$. Since every vertex in $C$ has at least one neighbor in $G'$ and at least two neighbors in $D$, every vertex in $C$ has degree at least $3$ in $G$ and if $C$ contains a vertex of degree $3$, then this is a light vertex. All vertices in $S_2^2$ are light, so to estimate $i(C)$ when $C$ is an exceptional heavy component, only the edge incidence count for $L_C$ will change. It remains true that at most two vertices of $G$ have degree $1$ in $C$, since $C$ is a cycle or a path, and now at most 17 vertices in $L_C$ are not incident to an edge of $M_0$. So 
\[ \sum_{u \in L_C} \left( |N(u) \cap D| + \frac{d_{G'}(u)}{2} \right) \geq \sum_{u \in L_C} \left( 2 + \frac{d_{C}(u) + 1}{2} \right) - \frac{17}{2} \geq \frac{7}{2}|L_C| - \frac{19}{2}.\]
Thus,
\[i(C) \geq \frac{7}{2}|L_C| + 4|H_C| - \frac{19}{2} \geq \frac{7}{2}|C| - \frac{17}{2}.\]

\item[Case 2:] $|V(C)| \geq 2$ and $C \in L$.

As in Case 1, if $C$ is a non-exceptional component, then we have
\[i(C) \geq \sum_{u \in L_C} \left(2 + \frac{d_C(u) + 1}{2}\right) + \sum_{c \in H_C} \left(3 + \frac{d_{C}(v)}{2} \right).\]
Now, either $|H_C| = 1$ or $H_C$ is empty. Note that since $2 + \frac{d_C(v) + 1}{2} < 3 + \frac{d_C(v)}{2}$, we have
\[i(C) \geq \sum_{u \in V(C)} 2 + \frac{d_C(u) + 1}{2}\]
regardless of the size of $H_C$. Now, if $C$ is a cycle, we have
\[i(C) \geq \sum_{u \in V(C)} 2 + \frac{2 + 1}{2} = \frac{7}{2}|C|.\]
If $C$ is a path, then we have $i(C) \geq \frac{7}{2}|C| - 1$.

If $C$ is an exceptional component then, as in Case 1, at most 17 vertices of $C$ either have degree $3$ or are in $S_2^2$ but not incident to an edge from $M_0$. We can then bound
\[i(C) \geq \sum_{u \in V(C)} \left(2 + \frac{d_C(u) + 1}{2}\right) - \frac{17}{2} \geq \frac{7}{2}|C| - \frac{19}{2}.\]
\end{enumerate}

Again, the goal is to show that very few components $C$ have $i(C) < \frac{7}{2}|C|$. Among non-exceptional components, only light paths fail to meet the desired value of $i(C)$. Our final step is thus to demonstrate that there are few non-exceptional light paths in $G'_0$. We first require a more careful description of heavy vertices in light paths.

Recall that if $C$ is a light path, then either $C$ contains no heavy vertices, or precisely one heavy vertex, say $v$. Note that if $v$ is a heavy vertex of $C$ which is incident to an edge of $M_0$, then we in fact have
\[i(C) = \sum_{u \in L_C} \left(2 + \frac{d_C(u) + 1}{2}\right) + 3 + \frac{d_C(v) + 1}{2} = \frac{5}{2}|C| + 1 + \sum_{u \in V(C)} \frac{d_C(u)}{2} = \frac{7}{2}|C|.\]
Thus, any light path $C$ with $i(C) < \frac{7}{2}|C|$ either contains no heavy vertex $v$, or its heavy vertex $v$ is not incident to an edge from $M_0$. Note that this implies $d_C(v) = d_{G'}(v)$, and in particular, $v$ is not adjacent to any heavy vertex in $V'$. Similarly, if $C$ contains a heavy vertex with more than $3$ neighbors in $D$, then $i(C) \geq \frac{7}{2}|C|$, so any light path $C$ with $i(C) < \frac{7}{2}|C|$ either contains no heavy vertex $v$, or its heavy vertex $v$ is adjacent to exactly $3$ vertices in $D$.

Now, recall that $D$ is precisely $V(G) \setminus N(x,y)$, so there is no triangle in $D$ containing $xy$. In particular, suppose $v$ is a heavy vertex in $G'$ which has no heavy neighbor in $G'$, such that $N(v) \cap D = \{x,y,z\}$. Then at most one neighbor of $v$ is adjacent to $z$, so $vz$ is an unrestricted edge. Among light paths, we will thus be able in large part to argue without distinguishing their heavy vertices.

Suppose $G_0'$ contains two non-exceptional light paths $P_1 = v_1,v_2,\dots, v_k$ and $P_2 = w_1, w_2, \dots, w_{\ell}$ such that $M_0$ contains no edges between $\{v_1, v_2,v_{k-1},v_k\}$ and $\{w_1,w_2,w_{\ell-1},w_{\ell}\}$. Note that we can find two such non-exceptional light paths if $G'_0$ contains more than 5 non-exceptional light paths in total. We will contradict that $G$ is properly rainbow $K_4$-saturated by showing that an edge between an endpoint of $P_1$ and an endpoint of $P_2$ may be added without creating any rainbow $K_4$-copy. 

To identify the correct place in which to add this edge, we will fix an orientation of both paths, as follows. Observe, since $c(v_1v_2) \neq 0$ and the $K_4$-copy on $\{x,y,v_1,v_2\}$ is not rainbow, either $c(v_1x) = c(v_2y)$ or $c(v_1y) = c(v_2x)$. We will say edge $v_iv_{i+1}$ of $P_1$ is \textit{left-oriented} if $c(v_ix) = c(v_{i+1}y)$, and \textit{right oriented} if not. Observe, if $v_1v_2$ is left-oriented, then $c(v_2y)$ cannot equal $c(v_3x)$, since $x$ is already adjacent to $v_1$ via an edge colored with $c(v_2y)$. Thus, if $v_1v_2$ is left-oriented, then $v_2v_3$ must also be left-oriented.
Inductively, if $v_1v_2$ is left-oriented, then in fact $v_{i}v_{i+1}$ must be left-oriented for all $i \leq k-1$, and similarly, if $v_1v_2$ is right-oriented, then $v_iv_{i+1}$ is right-oriented for every $i \leq k-1$.

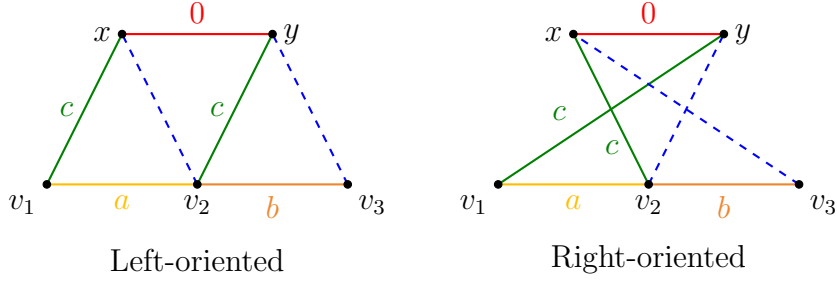
\begin{figure}[h]
\begin{center}
\begin{tikzpicture}

\draw[thick,red] (0,2) -- (2,2) node[pos=0.5, above]{0};

\draw[thick,amber] (-1,0) -- (1,0) node[pos=0.5, below]{$a$};
\draw[thick,cadmium] (1,0) -- (3,0) node[pos=0.5, below]{$b$};

\draw[thick, forest] (-1,0) -- (0,2) node[pos = 0.5, left]{$c$};

\draw[thick, forest] (1,0) -- (2,2) node[pos = 0.5, left]{$c$};

\draw[thick, dashed, blue] (1,0) -- (0,2);
\draw[thick, dashed, blue] (3,0) -- (2,2);

\filldraw (0,2) circle (0.05 cm);
\draw (0,2) node[left]{$x$};
\filldraw (2,2) circle (0.05 cm);
\draw (2,2) node[right]{$y$};

\filldraw (-1,0) circle (0.05 cm) node[below left]{$v_1$};
\filldraw (1,0) circle (0.05 cm) node[below]{$v_2$};
\filldraw (3,0) circle (0.05 cm) node[below right]{$v_3$};

\draw (1,-1) node{Left-oriented};

\draw[thick,red] (6,2) -- (8,2) node[pos=0.5, above]{0};

\draw[thick,amber] (5,0) -- (7,0) node[pos=0.5, below]{$a$};
\draw[thick,cadmium] (7,0) -- (9,0) node[pos=0.5, below]{$b$};

\draw[thick, forest] (5,0) -- (8,2) node[pos = 0.35, above left]{$c$};

\draw[thick, forest] (7,0) -- (6,2) node[pos = 0.25, left]{$c$};

\draw[thick, dashed, blue] (7,0) -- (8,2);
\draw[thick, dashed, blue] (9,0) -- (6,2);

\filldraw (6,2) circle (0.05 cm);
\draw (6,2) node[left]{$x$};
\filldraw (8,2) circle (0.05 cm);
\draw (8,2) node[right]{$y$};

\filldraw (5,0) circle (0.05 cm) node[below left]{$v_1$};
\filldraw (7,0) circle (0.05 cm) node[below]{$v_2$};
\filldraw (9,0) circle (0.05 cm) node[below right]{$v_3$};

\draw (7,-1) node{Right-oriented};

\end{tikzpicture}

\caption{Left and right orientations of $v_1v_2$. Note that in each case, to avoid a rainbow $K_4$-copy, the orientation of $v_2v_3$ must match that of $v_1v_2$.}
\end{center}
\end{figure}

Say $P_1$ is left-oriented if all of its edges are left-oriented, and right-oriented if all of its edges are right-oriented. Observe that orientation is purely a function of a path's vertex labelling: if $P_1$ is right-oriented, we may relabel its vertices, changing $v_i$ to $v_{n-i+1}$, to view $P_1$ as left-oriented. So, we relabel $P_1$ and $P_2$ if necessary to ensure that both are left-oriented. (Note that under this relabelling, it remains the case that no $M_0$ edge connects $\{v_1,v_2,v_{k-1},v_k\}$ and $\{w_1,w_2,w_{\ell - 1}, w_{\ell}\}$.) Now, we will add either $v_kw_1$ or $w_{\ell}v_1$ to $G$, choosing whichever edge contains fewer heavy vertices. Since $P_1,P_2$ are light, one of $v_kw_1$, $w_{\ell}v_1$ contains at most one heavy vertex. Without loss of generality, say we add $v_kw_1$.

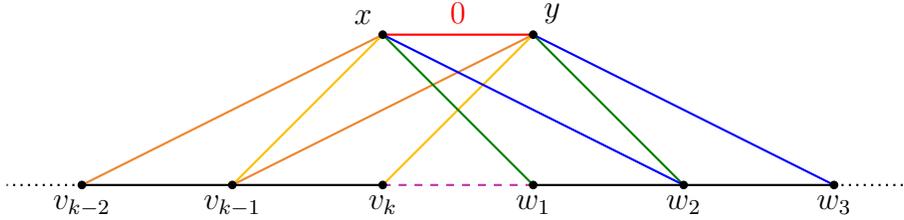
\begin{figure}[h]
\begin{center}
\begin{tikzpicture}

\draw[thick, cadmium] (3,2) -- (-1,0);
\draw[thick, cadmium] (5,2) -- (1,0);

\draw[thick, amber] (3,2) -- (1,0);
\draw[thick, amber] (5,2) -- (3,0);

\draw[thick, forest] (3,2) -- (5,0);
\draw[thick, forest] (5,2) -- (7,0);

\draw[thick, blue] (3,2) -- (7,0);
\draw[thick, blue] (5,2) -- (9,0);

\draw[thick, dashed, byzantine] (3,0) -- (5,0);

\draw[thick] (-1,0) -- (1,0);
\draw[thick] (1,0) -- (3,0);

\draw[thick, dotted] (-1,0) -- (-2,0);

\filldraw (-1,0) circle (0.05 cm) node[below]{$v_{k-2}$};
\filldraw (1,0) circle (0.05 cm) node[below]{$v_{k-1}$};
\filldraw (3,0) circle (0.05 cm) node[below]{$v_k$};

\draw[thick,red] (3,2) -- (5,2) node[pos=0.5, above]{0};

\draw[thick] (5,0) -- (7,0);
\draw[thick] (7,0) -- (9,0);

\draw[thick, dotted] (9,0) -- (10,0);

\filldraw (3,2) circle (0.05 cm);
\draw (3,2) node[above left]{$x$};
\filldraw (5,2) circle (0.05 cm);
\draw (5,2) node[above right]{$y$};

\filldraw (5,0) circle (0.05 cm) node[below]{$w_1$};
\filldraw (7,0) circle (0.05 cm) node[below]{$w_2$};
\filldraw (9,0) circle (0.05 cm) node[below]{$w_3$};

\end{tikzpicture}

\caption{$P_1$ and $P_2$, both left-oriented. We add the dashed edge $v_kw_1$ to $G$.}
\end{center}
\end{figure}

Because at most one of $v_k,w_1$ is heavy, the addition of $v_kw_1$ does not create any $K_4$-copy which intersects $D \setminus \{x,y\}$. In particular, if $P_1,P_2$ contain heavy vertices  $v_i, w_j$, say adjacent to $z_i,z_j \in D$ (possibly with $z_i = z_j$), then $c(v_iz_i)$ and $c(w_jz_j)$ are unrestricted, and even in $G + v_kw_1$, the edges $v_iz_i$ and $w_jz_j$ are not contained in any $K_4$-copy. We choose distinct new colors for each of $v_iz_i, w_jz_j$ to ensure that these edges do not share a color with any edge of $G$; we will also never reuse these new colors in subsequent recoloring steps. Thus, any recoloring of edges incident to $v_i,w_j$ will not conflict with $c(v_iz_i)$ or $c(w_jz_j)$. For the remainder of the proof, we may omit consideration of $v_iz_i$ and $w_jz_j$, if they exist.

The addition of $v_kw_1$ creates one $K_4$-copy, on $\{x,y,v_k, w_1\}$. Since no edge of $M_0$ connects $\{v_1,v_2,v_{k-1},v_k\}$ with $\{w_1,w_2,w_{\ell-1},w_{\ell}\}$, the addition of $v_kw_1$ does not create any triangle in $G'$, so in fact, $\{x,y,v_k, w_1\}$ span the unique $K_4$-copy containing $v_kw_1$. If $xv_k$ and $yw_1$ are both unrestricted edges, we can recolor so that $c(xv_k) = c(yw_1)$, using an edge-color not yet appearing in $G$. Coloring $v_kw_1$ with any legal color then creates a proper edge-coloring of $G + v_kw_1$ which is rainbow $K_4$-free, a contradiction.

Thus, suppose $xv_k$ is restricted. Since $P_1$ is left-oriented, the $K_4$-copy on $\{x,y, v_{k-1}, v_k\}$ has $c(xv_{k-1}) = c(xv_k)$, so for $xv_k$ to be restricted, $xv_k$ must be contained in another $K_4$-copy. This is only possible if $v_{k-2}v_{k} \in M_0$, a situation we depict in Figure~\ref{vk triangle}. Note that the coloring depicted in Figure~\ref{vk triangle} is without loss of generality; in particular, since $P_1$ is left-oriented, we have $c(xv_{k-2}) = c(xv_{k-1})$ and $c(xv_{k-1}) = c(yv_{k})$, none of which can be in $\{c(v_{k-2}v_{k-1}), c(v_{k-1}v_k)\}$ since the edge-coloring is proper.

\begin{figure}[h]

\begin{center}

\begin{tikzpicture}

\draw[thick, red] (0,2) -- (2,2) node[pos = 0.5, above]{\small 0};
\draw[thick, cadmium] (1,0) -- (3,0) node[pos = 0.5, above]{\small 1};
\draw[thick, amber] (-1,0) -- (1,0) node[pos = 0.5, above]{\small 2};
\draw[thick, red] (-1,0) to[bend right = 50] (3,0);
\draw[red] (1,-1) node[below]{\small 0};
\draw[thick, forest] (3,0) -- (2,2) node[pos = 0.5, right]{\small 3};
\draw[thick, forest] (1,0) -- (0,2) node[pos = 0.5, right]{\small 3};
\draw[thick, blue] (1,0) -- (2,2) node[pos = 0.5, left]{\small 4};
\draw[thick, blue] (-1,0) -- (0,2) node[pos = 0.5, left]{\small 4};

\filldraw (0,2) circle (0.05 cm) node[above]{$x$};
\filldraw (2,2) circle (0.05 cm) node[above]{$y$};

\filldraw (-1,0) circle (0.05 cm) node[below left]{$v_{k-2}$};
\filldraw (1,0) circle (0.05 cm) node[below]{$v_{k-1}$};
\filldraw (3,0) circle (0.05 cm) node[below right]{$v_k$};

\end{tikzpicture}

\caption{The structure of $G$ near $v_k$ if $xv_k$ is restricted}\label{vk triangle}

\end{center}

\end{figure}
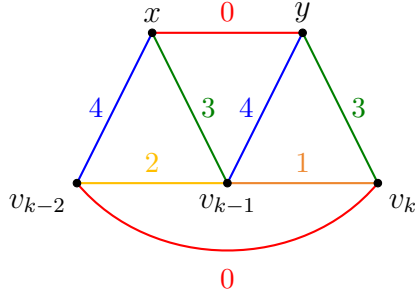

Observe that if $xv_k$ is restricted, then to avoid a rainbow $K_4$-copy, we must have $c(yv_{k-2}) = c(v_{k-1}v_k)$ (color $1$ in Figure~\ref{vk triangle}) and $c(xv_{k}) = c(v_{k-1}v_{k-2})$ (color $2$ in Figure~\ref{vk triangle}). Now, if $v_{k-2}v_{k-1}$ is only in the $K_4$-copies on $\{v_{k-2}, v_{k-1}, x,y\}$, $\{v_{k-2},v_{k-1},v_k,x\}$, and $\{v_{k-2},v_{k-1},v_k,y\}$, then we can recolor $v_{k-2}v_{k-1}$ and $xv_k$ simultaneously, so long as we maintain $c(v_{k-2}v_{k-1}) = c(xv_k)$. If $v_{k-2}v_{k-1}$ is in another $K_4$ copy in $G$, then this copy must include a triangle in $G'$, implying $v_{k-3}v_{k-1} \in M_0$. We depict this in Figure~\ref{vk component}.

\begin{figure}[h]

\begin{center}

\begin{tikzpicture}

\draw[thick, red] (-1,2) -- (1,2) node[pos = 0.5, above]{\small 0};
\draw[thick, cadmium] (1,0) -- (3,0) node[pos = 0.5, above]{\small 1};
\draw[thick, amber] (-1,0) -- (1,0) node[pos = 0.5, above]{\small 2};
\draw[thick, red] (-1,0) to[bend right = 50] (3,0);
\draw[red] (1,-1) node[below]{\small 0};
\draw[thick, red] (-3,0) to[bend right = 50] (1,0);
\draw[red] (-1,-1) node[below]{\small 0};
\draw[thick, forest] (3,0) -- (1,2) node[pos = 0.5, right]{\small 3};
\draw[thick, forest] (1,0) -- (-1,2) node[pos = 0.5, right]{\small 3};
\draw[thick, blue] (1,0) -- (1,2) node[pos = 0.5, left]{\small 4};
\draw[thick, blue] (-1,0) -- (-1,2) node[pos = 0.5, left]{\small 4};

\draw[thick, amber] (3,0) arc (0:127:2.5);
\draw[amber] (2.8,1.5) node{\small 2};
\draw[thick, cadmium] (-1,0) -- (1,2) node[pos=0.5,left]{\small 1};

\draw (-3,0) arc (180:53:2.5);
\draw (-3,0) -- (-1,2);
\draw (-3,0) -- (-1,0) node[pos = 0.5, above]{\small $b$};
\draw(-2.8,1.5) node{\small $a$};

\filldraw (-1,2) circle (0.05 cm) node[above]{$x$};
\filldraw (1,2) circle (0.05 cm) node[above]{$y$};

\filldraw (-3,0) circle (0.05 cm) node[below left]{$v_{k-3}$};
\filldraw (-1,0) circle (0.05 cm) node[below left]{$v_{k-2}$};
\filldraw (1,0) circle (0.05 cm) node[below right]{$v_{k-1}$};
\filldraw (3,0) circle (0.05 cm) node[below right]{$v_k$};

\end{tikzpicture}

\caption{The structure of $G$ near $v_k$ if both $v_{k-2}v_k$ and $v_{k-3}v_{k-1}$ are in $M_0$}\label{vk component}

\end{center}

\end{figure}
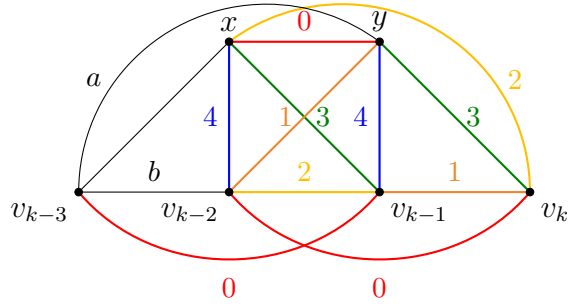

Now, consider $a = c(yv_{k-3})$ and $b = c(v_{k-3}v_{k-2})$. To maintain a proper coloring, $a \not\in \{0,1,3, 4\}$ and $b \not\in \{0,1,2,4\}$. So, since $y,v_{k-1}, v_{k-2}, v_{k-3}$ cannot span a rainbow $K_4$-copy in $G$, we must have $a = 2$. Because $P_1$ is left-oriented, we know that $a = c(yv_{k-3})$ is equal to $c(xv_{k-4})$, if $v_{k-4}$ exists. However, such a color assignment would be impossible if $a = 2$, since $x$ is already incident to $v_k$ via an edge of color $2$. 

Thus, if $v_{k-3}v_{k-1}$ is in $M_0$, we conclude that $P_1$ is a four-vertex light path with $v_{k-3}v_{k-1}$ and $v_{k-1}v_k$ in $M_0$. In this case, $P_1$ is not only a component of $G'_0$ but of $G'$ and, because $P_1$ is light, every $K_4$-copy in $G + v_{k-3}v_k$ which intersects $P_1$ must use only vertices from $V(P_1) \cup \{x,y\}$.
We can thus add the edge $v_{k-3}v_k$, since $P_2 \vee K_4$ can be properly edge-colored without creating a rainbow $K_4$-copy, a contradiction.

Thus, either $xv_k$ is unrestricted, or $c(xv_k) = c(v_{k-2}v_{k-1})$ and we can freely change the value of $c(xv_k)$ so long as we also change $c(v_{k-2}v_{k-1})$ to maintain $c(xv_k) = c(v_{k-2}v_{k-1})$. Analogously, either $yw_1$ is unrestricted, or $c(yw_1) = c(w_2w_3)$ and we can freely change the value of $c(yw_1)$ so long as we also change the value of $c(w_2w_3)$. Also note that because $P_1,P_2$ are separate components, the edges $xv_k, v_{k-1}v_{k-2}, yw_1,$ and $w_2w_3$ are pairwise vertex disjoint, so it is possible to maintain a legal edge-coloring if some or all of their colors are equal. We recolor so that $c(xv_k) = c(yw_1)$ via a color not yet used in $G$, also recoloring $v_{k-2}v_{k-1}$ and $w_2w_3$ if necessary. Under this new edge-coloring, $G$ remains rainbow $K_4$-free. Moreover, we have constructed this new coloring to allow the addition of $v_kw_1$ in any legal color without creating a rainbow $K_4$-copy.

Thus, if we can find $P_1,P_2$ in $G'_0$, then $G$ is not properly rainbow $K_4$-saturated, a contradiction. If we cannot find $P_1,P_2$ as above, then $G'_0$ contains at most $5$ non-exceptional light paths. We know that $G'_0$ contains at most 17 exceptional components in total, so contains at most 22 components $C$ with $i(C) < \frac{7}{2}|C|$. We have seen that every component has $i(C) \geq \frac{7}{2}|C| - \frac{19}{2} > \frac{7}{2}|C| - 10$, so in total
\[e(G) \geq \sum_{C} i(C) \geq \sum_{C} \frac{7}{2} |C| - 22\cdot 10 = \frac{7}{2}(n - |D|) - 220 \geq \frac{7}{2}(n - 2\alpha n) - 220 > \frac{7}{2}n - 8\alpha n\]
since we have $\alpha n > 220$.
\end{proof}

\section{Paths}\label{sec:paths}

In this section, we prove Theorem~\ref{paths thm}. We begin by proving the desired lower bound in the following proposition. We remark that this argument can be adapted to show that $\mathrm{sat^*}(n,T) \geq n-1$ for all trees $T$ with diameter at least $4$.

\begin{prop}\label{lem:P5-lowerbound-acycliccomps}
    If a graph $G$ is rainbow $P_k$-saturated for $k \geq 5$, then $G$ has at most one acyclic component. Consequently, $\sat^*(n,P_k) \geq n-1$.
\end{prop}
\begin{proof}
    We first show that $G$ does not contain two acyclic components of diameter at most $2$. For the sake of a contradiction, assume two such components exist. Note that these two components are stars, and adding an edge between their centers will create no new $P_k$-copies. Thus, $G$ is not rainbow $P_k$-saturated, a contradiction. It follows that $G$ has at most one acyclic component of diameter at most $2$. If this component has at least $3$ vertices, then adding an edge between two leaves will not create any new $P_k$-copies. So again, $G$ is not rainbow $P_k$-saturated, a contradiction. Thus, such a component has at most $2$ vertices, meaning it is isomorphic to $K_1$ or $K_2$.
    
    We now show that, in fact, $G$ cannot contain two acyclic components of any diameter. Suppose towards a contradiction that $G$ has two acyclic components, $T_1$ and $T_2$. By the above observations, one of $T_1$ or $T_2$ must have diameter at least $3$. Without loss of generality, say this is $T_1$. 
    Fix $u\in V(T_1)$ and $x\in V(T_2)$ to be endpoints of longest shortest paths (i.e., paths which realize the diameters of $T_1$ and $T_2$). Then $u$ and $x$ are necessarily leaves. 
    By the above, we may assume $T_1$ is a tree of diameter at least 3, while $T_2$ may be either $K_1$, $K_2$, or a tree of diameter at least 3. 
    
    Let $v$ be the unique neighbor of $u$. Since $T_1$ has diameter at least $3$, $v$ must have exactly one non-leaf neighbor, say $w$. Let $c$ be a rainbow $P_k$-free proper edge-coloring of $G$. Without loss of generality, $c(uv) = 1$ and $c(vw) = 2$. Now, we consider the edge-coloring of $T_2$. If $T_2$ is not a $K_1$-copy, let $y$ be the neighbor of $x$, and if $T_2$ is not a $K_2$-copy, let $z$ be the unique non-leaf neighbor of $y$. By relabeling edge-colors within $T_2$ only, we can set $c(xy) = 1$ (if $y$ exists) and $c(yz) = 2$ (if $z$ exists).

    \begin{figure}[h]

\begin{center}

\begin{tikzpicture}

\draw[thick, red] (2,0) -- (4,0) node[pos = 0.5, above]{1};
\draw[thick, blue] (0,0) -- (2,0) node[pos = 0.5, above]{2};

\draw[thick, red] (6,0) -- (8,0) node[pos = 0.5, above]{$a$};
\draw[thick, blue] (8,0) -- (10,0) node[pos = 0.5, above]{$b$};

\draw[thick, blue, dashed] (4,0) -- (6,0) node[pos = 0.5, above]{2};

\draw[thick, forest] (2,0) -- (1.56,1);
\draw[thick, amber] (2,0) -- (2,1);
\draw[thick, cadmium] (2,0) -- (2.44,1);

\draw[thick, fill = blue!20!white] (2,1.25) circle (0.5 cm) node{\tiny Leaves};

\draw[thick, byzantine] (8,0) -- (7.56,1);
\draw[thick, forest] (8,0) -- (8,1);
\draw[thick, cyan] (8,0) -- (8.44,1);

\draw[thick, fill = blue!20!white] (8,1.25) circle (0.5 cm) node{\tiny Leaves};

\draw (-1,0) -- (0,0);

\draw (10,0) -- (11,0);

\filldraw (0,0) circle (0.05 cm) node[below]{$w$};
\filldraw (2,0) circle (0.05 cm) node[below]{$v$};
\filldraw (4,0) circle (0.05 cm) node[below]{$u$};

\filldraw (6,0) circle (0.05 cm) node[below]{$x$};
\filldraw (8,0) circle (0.05 cm) node[below]{$y$};
\filldraw (10,0) circle (0.05 cm) node[below]{$z$};

\end{tikzpicture}

\caption{Leaves in $T_1$ and $T_2$ which can be connected by an edge of color $2$. If necessary, we can relabel colors within $T_2$ only so that $a = 1$ and $b = 2$}\label{tree components}

\end{center}

\end{figure}
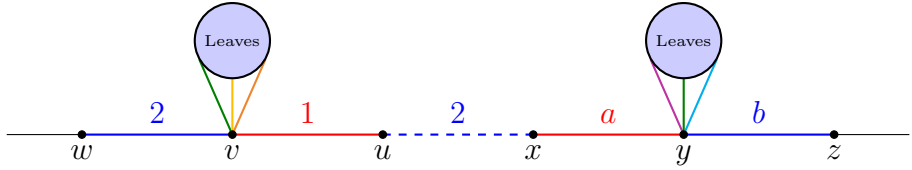

    Now we claim we can add the edge $ux$ with $c(ux) = 2$ without creating any rainbow $P_k$-copy. Indeed, observe that since $k \geq 5$, any copy of $P_k$ which contains $ux$ either contains one of $vw$ or $yz$, or else both of $uv$ and $xy$. Thus, $G + ux$ can be properly edge-colored without rainbow $P_k$-copies, contradicting the assumption that $G$ is properly rainbow $P_k$-saturated.
\end{proof}

We now work to show the upper bound in Theorem~\ref{paths thm}. The substance of this argument will be showing that a particular construction is properly rainbow $P_k$-saturated. In order to do this, we must define a family of graphs upon which our construction is based, and develop an understanding of the proper edge-colorings of this family.

\begin{defn}
For $w\geq 3$, let $V = \{0,1\}^w/\equiv,$ where $\equiv$ is the equivalence relation $u\equiv v$ iff $u_i = 1-v_i$ for each $i\in [w]$. Let $uv\in E$ if and only if $u,v\in V$ have respective representatives that differ in exactly one bit. Define $H_w/2 := (V,E)$. 
\end{defn}

We note that $H_w/2$ is precisely the $(w-1)$-dimensional hypercube $H_{w-1}$ with edges joining antipodal vertices, sometimes called the \textit{folded hypercube}. Moreover, $H_w/2$ is $w$-regular, and $\operatorname{diam}(H_w/2) = \floor{w/2}$. The folded hypercube has previously been examined in the study of rainbow Tur\'an numbers for paths; in the cases where the rainbow Tur\'an number of $P_k$ is known, folded hypercubes provide an extremal construction.

The following lemma describes a construction of Johnston and Rombach \cite{Johnston2019LowerBF} for a coloring of $H_w/2$ which avoids rainbow $P_{w+1}$-copies. In this lemma and the rest of this section, we will find it convenient to describe the edges of $H_w/2$ by the component in which their endpoints differ. We refer to the elements of the binary strings defining the vertices as \emph{bits} and say that edge $e\in E(H_w/2)$ corresponds to an $i^{th}$ \emph{bit-flip} if (some representations of) its endpoints differ in precisely bit $i$. 

 \begin{lem}\label{lem:H_w coloring exists}
 For $w\geq 3$, $H_w/2$ has a proper edge-coloring that is rainbow $P_{w+1}$-free.
 \end{lem}
 \begin{proof}
Define a coloring $c: E\to \mathbb{N}$ by setting $c(e) = i$ if $e$ corresponds to an $i^{th}$ bit-flip. Then for any $u\in V$, a rainbow walk starting from $u$ with $w$ edges flips every bit exactly once, and therefore forms a cycle. Hence, $H_w/2$ has no rainbow $P_{w+1}$ under this coloring.
 \end{proof}

While we have exhibited a proper edge-coloring of $H_w/2$ which is rainbow $P_{w+1}$-free, this is not necessarily the unique such coloring. In order to use the family of folded hypercubes to construct the desired properly rainbow $P_k$-saturated family of graphs, we must understand the set of all rainbow $P_{w+1}$-free proper edge-colorings of $H_w/2$. Towards this end, we show that every rainbow $P_{w+1}$-free proper edge-coloring of $H_w/2$ contains a rainbow cycle of particular structure.

\begin{defn}
A total bit-flip cycle (TBFC) is a cycle of length $w$ in $H_w/2$ in which every edge corresponds to a different bit-flip. A total bit-flip rainbow cycle (TBFRC) is a TBFC with a rainbow edge-coloring.
\end{defn}

\begin{lem}\label{lem:H_w has TBFRC}
For $w\geq 3$, every proper edge-coloring of $H_w/2$ that is rainbow $P_{w+1}$-free contains a TBFRC.
\end{lem}
\begin{proof}
We construct a rainbow path $P_k = (x_0,\ldots, x_k)$ for all $0\leq k\leq w-1$ iteratively. Set $P_0 = (x_0)$, for some fixed vertex $x_0$. Given $P_k$, the manner in which we extend to $P_{k+1}$ will depend upon the value of $k$. We will say that a color $c$ is \textit{new} for $P_k$ if $c \neq c(x_{i-1}x_{i})$ for any $i \in \{1, \dots, k\}$. Let $d:= \operatorname{diam}(H_w/2) = \floor{w/2}$. While $k \leq d$, we say a color $c$ is \textit{bad} for $P_k$ if $c$ satisfies the following conditions:

\begin{enumerate}
\item $c = c(x_{i-1}x_{i})$ for some $i \in \{1, \dots, k\}$;
\item For some vertex $y \not\in \{x_0, \dots, x_{k-1}\}$, $x_ky$ is an edge with $c(x_ky) = c$.
\end{enumerate}

For $0\leq k\leq d$, extend $P_k$ to $P_{k+1}$ by adding an edge $x_k x_{k+1}$ such that $c(x_k x_{k+1})$ is new (for convenience we may assume this is color $k+1$) and $x_k x_{k+1}$ corresponds to a yet unused bit-flip. Such an edge exists, since the number of bad colors at this step is at most $k-1$, and the edges of $P_k$ correspond to $k$ bit-flips, leaving $w-k$ new bits available to flip. Indeed, if $k \leq d$, then $2k-1 \leq 2\floor{w/2}-1$, and thus, we have $k-1 < w-k$. 

In this manner, we produce a rainbow path $P_{d+1}$ of length $d+1$, where each edge uses a different bit-flip. Now, for $k \geq d + 1$, we continue an extension process, obtaining $P_{k+1}$ by adding an edge $x_kx_{k+1}$ to $P_k$ so that $c(x_kx_{k+1})$ is new and $x_kx_{k+1}$ corresponds to a specified type of bit-flip. However, to ensure that such an edge $x_kx_{k+1}$ exists, we must somewhat relax the condition on which bit-flips may correspond to $x_kx_{k+1}$. To describe this new condition, we introduce the following notation. 

Given $d+1 \leq k \leq w-1$, suppose $P_k$ has been defined. We denote by $S_d$ the subpath $(x_0, \dots, x_{d})$ of $P_k$, and for $d+1 \leq j \leq k$, we denote by $S_j$ the subpath $(x_{2j - w}, \dots, x_j)$ of $P_k$. We say that bit-flip $i$ is \textit{used by} $S_j$ if some edge of $S_j$ corresponds to bit-flip $i$. Rather than extending $P_k$ by adding an edge $x_kx_{k+1}$ which does not correspond to the same bit-flip as any previous edge in $P_k$, we will only demand that $x_kx_{k+1}$ corresponds to a bit-flip not used by $S_k$. For $d \leq j < k$, we will denote by $S_j^+$ the concatenation of $S_j$ and $x_jx_{j+1}$. Thus, to determine an acceptable choice for $x_{j+1}$ as we are building our path, we ensure that $x_jx_{j+1}$ does not repeat a bit-flip used by $S_j$; once we add $x_{j+1}$ to our path, the subpath $S_j^+$ contains no repeated bit-flips.

\begin{claim}\label{good path}
For $d+1 \leq k \leq w-1$, there exists a rainbow path $P_k = (x_0, \dots x_k)$ in $H_w/2$ such that for every $j$ with $d \leq j \leq k-1$, $S_j^+$ contains no two edges corresponding to the same bit-flip. Moreover, $x_ix_j \not \in E(H_w/2)$ if $i+1 \not\equiv j \mod w$. 
\end{claim}

\begin{proof}[Proof of Claim~\ref{good path}]

We proceed by induction. For $k = d+1$, we take $P_{d+1}$ as described above. By construction, $P_{d+1}$ is rainbow, and each edge corresponds to a different bit-flip. It is immediate that
$S_d^+$ contains no two edges corresponding to the same bit-flip, and that two vertices $x_i,x_j$ either are adjacent in $P_{d+1}$ or
are connected by a path corresponding to $\ell \geq 2$ distinct bit-flips.
Note that
if $2 \leq \ell < w-1$, then $x_i,x_j$
are not adjacent in $H_w/2$. If 
$\ell = w-1$, then $x_i = x_0$ and $x_j = x_{d+1}$,
and we must have $d+2 \equiv 0 \mod w$. Thus, the statement holds for $k = d +1$. 

Let $d+1 < k \leq w-1$ and suppose that there exists a rainbow path $P_{k-1} = (x_0, \dots, x_{k-1})$ satisfying all conditions in the statement of Claim~\ref{good path}. We shall obtain a rainbow path $P_k$ by extending $P_{k-1}$. Let $N$ be the set of neighbors of $x_{k-1}$ which are not on $P_{k-1}$. Since $k \leq w-1$, by the inductive hypothesis, $x_{k-1}$ is adjacent to no vertex on $P_{k-1}$ except for $x_{k-2}$, so $|N| = w - 1$. Note that at most $k-2$ elements of $N$ are adjacent to $x_{k-1}$ via an edge whose color is not new for $P_{k-1}$ (since the color of $x_{k-2}x_{k-1}$ is not new for $P_{k-1}$, but is not used for any edge between $x_{k-1}$ and $N$). Moreover, exactly $w- k + 1$ bit-flips are used by $S_{k-1}$, so at most $w - k$ elements of $N$ are adjacent to $x_{k-1}$ via an edge corresponding to a bit-flip used by $S_{k-1}$ (since the bit-flip corresponding to $x_{k-2}x_{k-1}$ does not correspond to any edge from $x_{k-1}$ to $N$). In total, there are at most $(k-2) + (w - k) = w - 2$ elements $x$ of $N$ such that either $c(x_{k-1}x)$ is not new for $P_{k-1}$ or $x_{k-1}x$ corresponds to a bit-flip used by $S_{k-1}$. Thus, there exists an element $x_{k}$ of $N$ such that $c(x_{k-1}x_k)$ is new for $P_{k-1}$ and $x_{k-1}x_k$ does not correspond to a bit-flip used by $S_{k-1}$. We set $P_k = (x_0, \dots, x_{k-1},x_k)$.

Now, it remains to show that $P_k$ satisfies all conditions in Claim~\ref{good path}. By construction, $P_k$ is rainbow. For $j < k-1$, the inductive hypothesis on $P_{k-1}$ implies that $S_j^+$ contains no two edges corresponding to the same bit-flip. By construction, $S_{k-1}^+$ contains no two edges corresponding to the same bit-flip, 
since $x_kx_{k+1}$ does not correspond to a bit-flip used by any edge of $S_{k-1}$, and since $S_{k-1}$ is contained in $S_{k-2}^+$, implying that no two edges of $S_{k-1}$ correspond to the same bit-flip. 
Thus, we only need verify that $x_ix_j \not\in E(H_w/2)$ if $i+1 \not\equiv j \mod w$. Again, the inductive hypothesis on $P_{k-1}$ guarantees this if both $i,j$ are smaller than $k$. Thus, we consider pairs $x_i,x_k$. There are several cases to consider. 

\begin{enumerate}

\item[Case 1:] $i = 0$ and $k \leq w - 2$

Since $S_{d}^+$ has no repeated bit-flips, we have $d(x_0,x_{d+1}) = w-(d+1)$. Then since $d(x_0,x_{d+1}) \leq d(x_0,x_k) + d(x_k,x_{d+1})$, we have \[w-(d+1) \leq d(x_0,x_k) + k-(d+1).\] This implies $2\leq w-k\leq d(x_0,x_k)$. Therefore, $x_0x_k\not\in E(H_w/2)$.

\item[Case 2:] $i = 1$

By assumption, $S_{d+1}^+$ has no repeated bit-flips and has length $w-d$, so we have $d(x_{2d-w+2}, x_{d+2}) = d.$ Thus, $d(x_{2d-w+2},x_{d+2}) \leq d(x_{2d-w+2}, x_1) + d(x_1,x_k) + d(x_k,x_{d+2})$, which implies \[d\leq (2d-w+1) + d(x_1,x_k) + k-(d+2).\] As a result, we have $w-k+1\leq d(x_1,x_k)$, and hence $2\leq d(x_1,x_k)$. Therefore, $x_1x_k \not\in E(H_w/2)$.

\item[Case 3:] $2 \leq i \leq 2(k-1) - w$

Let $j = \floor{\frac{i+w}{2}}$. Note that $j \geq \floor{\frac{2+w}{2}} = d+1$, and $j \leq \floor{\frac{2(k-1)-w+w}{2}} = k-1$. Hence by assumption, $S_j^+$ has no repeated bit-flips, and is a path from $x_{2j-w}$ to $x_{j+1}$ with length $w-j+1$. We also note that $i-1\leq 2j-w\leq i$ and $i< j+1$. If $i=2$ and $w$ is odd, then $j = d + 1$ and $S_j^+$ contains $w - (d+1) + 1 = w-d = d + 1$ edges. In this case, $2j - w = 1$ and $i = 2$, so $d(x_i,x_{j+1}) = d$, and we note that $d = j+1 - i$. If $w$ is even or $i \geq 3$, then $S_j^+$ contains $w-j+1 \leq d$ edges, and it follows that $d(x_i,x_{j+1}) = j+1 - i$.
We combine these observations with the inequality $d(x_i,x_{j+1}) \leq d(x_i,x_k) + d(x_k,x_{j+1})$, which implies \[j+1-i \leq d(x_i,x_k) + (k-(j+1)).\] Thus, we have $2 + 2j - i \leq d(x_i, x_k) + (w-1)$, which gives $2\leq d(x_i,x_k).$ Therefore, $x_ix_k \not\in E(H_w/2)$.

\item[Case 4:] $2(k-1) - w < i \leq k -2$
By assumption, $S_{k-1}^+$ has no repeated bit-flips and is a path from $x_{2(k-1)-w}$ to $x_k$ with length $w-k+2$. Then $x_i$ lies on this path, with $d(x_i,x_k) = k-i\geq 2$. Therefore, $x_ix_k \not\in E(H_w/2)$.
\end{enumerate}
\end{proof}

With Claim~\ref{good path} established, we construct a TBFRC in $H_w/2$ as follows. Let $P_{w-1} = (x_0, \dots x_{w - 1})$ be a rainbow path satisfying the conditions of Claim~\ref{good path}. Without loss of generality, $c(x_{i-1}x_{i}) = i$ for each $i \in \{1, \dots, w - 1\}$. Since the only neighbors of $x_{w-1}$ on $P_{w-1}$ are $x_{w-2}$ and potentially $x_0$, there are at least $w-2$ ways to extend $P_{w-1}$ to a path with $w$ edges. By hypothesis, $H_w/2$ is rainbow $P_{w+1}$-free, so each of these $w-2$ edges must use colors from $\{1, \dots, w-2\}$. This also means that there must exist an edge $x_{w-1} x_0$, and it must use a new color, $w$. Let $C_w:= P_{w-1} + x_{w-1}x_0$. Thus, $C_w$ is a rainbow cycle of length $w$. It remains to check that $C_w$ is a TBFC. 

Fix a representation of $x_0$. For $1 \leq i \leq w$, if $x_{i-1}$ and $x_i$ differ in bit $j$, then fix the representation of $x_i$ obtained by flipping the $j$th bit of $x_{i-1}$. In this manner, we select representations for $x_0,\dots,x_w$ so that for each $1\leq i \leq j$, there is precisely one bit of difference between $x_{i-1}$ and $x_i$. With these representations fixed, we proceed as follows. Since $x_0x_{w-1} \in E(H_w/2)$, we have that (the fixed representations of) $x_0$ and $x_{w-1}$ either differ in exactly one bit or differ in $w-1$ bits. Observe that if $x_{w-1}$ differs from $x_0$ in $w-1$ bits, then $C_{w}$  is a TBFC. Thus, it suffices to show that $x_0$ and $x_{w-1}$ differ in more than one bit. Note that $x_0$ and $x_{w-1}$ differ in bit $i$ if and only if an odd number of edges in $P_{w-1}$ correspond to the $i$th bit-flip. So if $x_0$ and $x_{w-1}$ differ in exactly one bit, then the number of edges of $P_{w-1}$ corresponding to bits in which $x_0$ and $x_{w-1}$ do not differ is even, and the number of edges corresponding to the unique bit in which they do differ is odd. If $w$ is odd, then $P_{w-1}$ has an even number of edges, and it is thus impossible for $x_0$ and $x_{w-1}$ to differ in exactly one bit. If $w$ is even, consider the path $P = (x_d, x_{d+1}, \dots x_0)$ on $C_w$. Because no bit-flips are repeated within $S_d^+$, we know that $x_d$ in fact has distance $d = \frac{w}{2}$ from $x_0$. Since $P$ is a path of length $d$ from $x_0$ to $x_d$, $P$ must also contain no two edges which correspond to the same bit-flip. Thus, the bit-flip corresponding to $x_dx_{d+1}$ is repeated nowhere in $P_{w-1}$. Thus, if $x_0$ and $x_{w-1}$ differ only in one bit, then this is the same bit in which $x_{d}$ differs from $x_{d+1}$; moreover, every bit-flip which corresponds to another edge in $P_{w-1}$ in fact corresponds to at least two edges of $P_{w-1}$. 

Now, we partition $P_{w-1}$ into two sub-paths: $P_1 = (x_0, \dots , x_{d})$ and $P_2 = (x_{d+1}, \dots x_{w-1})$. Since both $P_1$ and $P_2$ contain no repeated bit-flips, it must be the case that every bit-flip corresponding to an edge of $P_1$ also corresponds to an edge of $P_2$. However, $P_1$ has length $d$ and $P_2$ has length $d-2$, and thus some bit-flip corresponding to an edge of $P_1$ does not correspond to an edge of $P_2$. We conclude that $x_0$ and $x_{w-1}$ differ in more than one bit; thus, $C_w$ is a TBFRC.
\end{proof}

Now, the existence of a TBFRC in any rainbow $P_{w+1}$-free edge-coloring of $H_w/2$ will imply that in fact, all rainbow $P_{w+1}$-free edge-colorings of $H_w/2$ are equivalent up to relabeling of colors. To see this, we also require the following lemma.

\begin{lem}\label{lem:TBFRC extends to H_w}
For $w\geq 5$, a TBFRC in $H_w/2$ can be uniquely extended to a rainbow $P_{w+1}$-free edge-coloring of $H_w/2$.
\end{lem}
\begin{proof} Suppose that we have a TBFRC $C_w := (x_1,x_2,\ldots,x_w,x_1)$, where $x_1$ is an arbitrary starting point. Reindexing if necessary, we may assume that $x_ix_{i+1}$ corresponds to the $i$th bit-flip in $H_w/2$. We identify any TBFC starting at $x_1$ with the permutation on [$w$] given by the order of bit-flips corresponding to the edges of this cycle; thus, $C_w$ corresponds to the identity permutation. We may also assume without loss of generality that $c(x_ix_{i+1})=i$ (with subscripts taken modulo $w + 1$). 

Note that if $x$ is on a TBFRC $C$, and $y\in N(x)$ is not on $C$, then $c(xy)$ must be a color appearing on $C$ to avoid a rainbow $P_w$-copy. Also note that the only vertices on $C$ which are adjacent to $x$ are the two neighbors of $x$ along $C$, since a vertex at distance at least $2$ from $x$ in $C$ differs from $x$ in at least two bit-flips. Thus, $x$ is incident to only edges receiving colors used on $C$. Since $d(x) = w$, this implies that $x$ is incident to edges of all colors used on $C$.

We now consider TBFC's in $H_w/2$ which interact with $C_w$. For a fixed $2 \leq i \leq w$, let $y$ be the vertex such that $x_{i-1}y$ corresponds to the 
$i$th bit-flip in $H_w/2$. Note that $yx_{i+1}$ is an edge of $H_w/2$ corresponding to the 
$(i-1)$st bit-flip. So, \[C(i, i-1) := (x_1, \dots, x_{i-1}, y, x_{i+1}, \dots x_w, x_1)\] is a TBFC not containing $x_i$, which corresponds to the transposition $(i,i-1)$.
Since we assume $w \geq 5$, note also that $N(x_{i-1}, x_{i+1})$ is exactly equal to $\{x_i, y\}$. (If $w = 4$, this statement does not hold; common neighborhoods in $H_w/2$ have size $3$.) We now consider the edge-coloring of $C(i,i-1)$. 

We claim that $c(x_{i-1}y) = i$ and $c(yx_{i+1}) = i-1$. Indeed, if this did not hold, then by the above observations on the edges incident to $x_{i-1}$ and $x_{i+1}$, there would exist edges $x_{i-1}u$ and $x_{i+1}v$ such that $u \neq v$ and $u,v \not\in V(C(i,i-1))$ with $c(x_{i-1}u) = i$ and $c(x_{i+1}v) = i-1$. (We allow for the possibility that one of $u,v$ is equal to $y$.) This would yield that $(u, x_{i-1}, x_{i-2}, \dots, x_1, x_w, \dots, x_{i+1}, v)$ is a rainbow $P_{w+1}$-copy in $H_w/2$, a contradiction. Thus, $c(x_{i-1}y) = i$ and $c(yx_{i+1}) = i-1$, so $C(i, i-1)$ is a TBFRC. Moreover, the edge-coloring of $C(i,i-1)$ is uniquely determined, with the property that each edge of $C(i,i-1)$ has the same color as the bit-flip to which it corresponds.

Repeatedly applying the above argument, if $C'$ is a TBFC containing $x_1$ which corresponds to a product of transpositions of the form $(i,i-1)$, then $C'$ has a rainbow edge-coloring determined by the permutation to which it corresponds. That is, if $u,v$ are adjacent on $C'$ and differ in bit $i$, then $c(u,v) = i$. Note that every TBFC containing $x_1$ corresponds to some permutation of $[w]$, and it is well-known that every permutation of $[w]$ can be written as the product of permutations of the form $(i,i-1)$. Thus, if $uv$ is an edge of $H_w/2$ such that $uv$ is in some TBFC containing $x_1$, then $c(uv)$ is equal to the bit in which $u,v$ differ. We observe that given $x_1$ and $uv$ (possibly with $x_1 = u$), it \textit{is} possible to construct a TBFC beginning at $x_1$ and containing $uv$. Indeed, let $P_u$ be a shortest path between $x_1$ and $u$, and $P_v$ a shortest path between $x_1$ and $v$. Thus, $P_u$ and $P_v$ do not repeat bit-flips. Without loss of generality, we may assume $|e(P_u)| \leq |e(P_v)|$. Suppose that $uv$ corresponds to bit-flip $i$, and there is also an edge of $P_u$ corresponding to bit-flip $i$. Let $B$ be the set of bit-flips different from $i$ which correspond to edges of $P_u$. Note that $v$ must differ from $x_1$ in precisely the bits of $B$, so there is a path of length $|B| < |e(P_u)|$ from $x_1$ to $v$, contradicting the assumption that $|e(P_u)| \leq |e(P_v)|$. Thus, $P_u + uv$ is a path with no repeated bit-flips containing $x_1$ and $uv$. It is clear that $P_u + uv$ can be extended to a TBFC. Thus, $uv$ is colored by the bit in which $u,v$ differ.
\end{proof}

We are now equipped to prove our main result, Theorem~\ref{thm: unique Pw+1-free coloring of Hw/2}, on the edge-colorings of $H_w/2$.

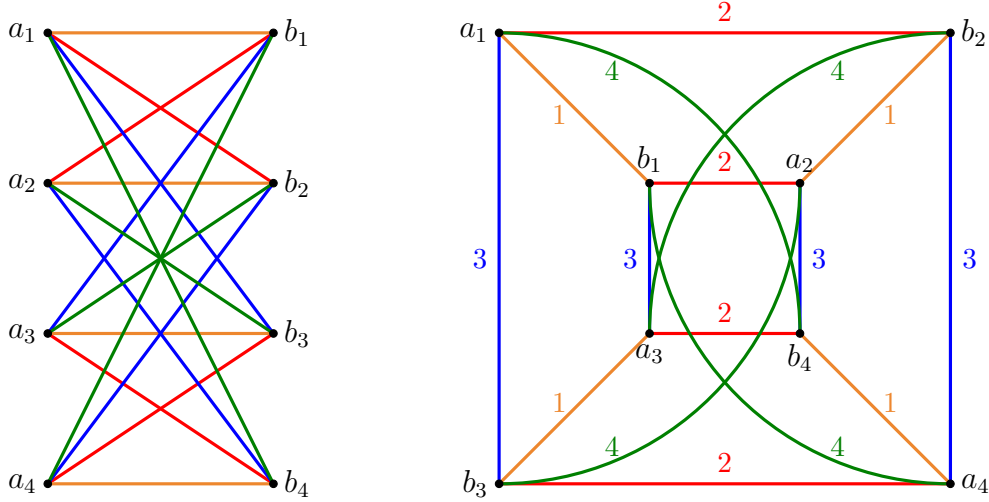
\begin{figure}[h!]

\begin{center}

\begin{tikzpicture}

\draw[very thick, red] (0,6) -- (3,4);
\draw[very thick, red] (0,4) -- (3,6);
\draw[very thick, red] (0,2) -- (3,0);
\draw[very thick, red] (0,0) -- (3,2);

\draw[very thick, cadmium] (0,6) -- (3,6);
\draw[very thick, cadmium] (0,4) -- (3,4);
\draw[very thick, cadmium] (0,2) -- (3,2);
\draw[very thick, cadmium] (0,0) -- (3,0);

\draw[very thick, blue] (0,6) -- (3,2);
\draw[very thick, blue] (0,2) -- (3,6);
\draw[very thick, blue] (0,4) -- (3,0);
\draw[very thick, blue] (0,0) -- (3,4);

\draw[very thick, forest] (0,6) -- (3,0);
\draw[very thick, forest] (0,0) -- (3,6);
\draw[very thick, forest] (0,4) -- (3,2);
\draw[very thick, forest] (0,2) -- (3,4);

\filldraw (0,0) circle (0.05 cm) node[left]{$a_4$};

\filldraw (0,2) circle (0.05 cm) node[left]{$a_3$};

\filldraw (0,4) circle (0.05 cm) node[left]{$a_2$};

\filldraw (0,6) circle (0.05 cm) node[left]{$a_1$};

\filldraw (3,0) circle (0.05 cm) node[right]{$b_4$};

\filldraw (3,2) circle (0.05 cm) node[right]{$b_3$};

\filldraw (3,4) circle (0.05 cm) node[right]{$b_2$};

\filldraw (3,6) circle (0.05 cm) node[right]{$b_1$};

\draw[very thick, cadmium] (6,0) -- (8,2) node[pos = 0.4, above]{\small 1};
\draw[very thick, cadmium] (12,0) -- (10,2) node[pos = 0.4, above]{\small 1};
\draw[very thick, cadmium] (6,6) -- (8,4)node[pos = 0.4, below]{\small 1};
\draw[very thick, cadmium] (12,6) -- (10,4)node[pos = 0.4, below]{\small 1};

\draw[very thick, red] (6,6) -- (12,6) node[pos = 0.5, above]{\small 2};
\draw[very thick, red] (8,4) -- (10,4) node[pos = 0.5, above]{\small 2}; 
\draw[very thick, red] (6,0) -- (12,0) node[pos = 0.5, above]{\small 2};
\draw[very thick, red] (8,2) -- (10,2) node[pos = 0.5, above]{\small 2}; 

\draw[very thick, blue] (6,0) -- (6,6)  node[pos = 0.5, left]{\small 3};
\draw[very thick, blue] (8,2) -- (8,4) node[pos = 0.5, left]{\small 3};
\draw[very thick, blue] (12,0) -- (12,6) node[pos = 0.5, right]{\small 3};
\draw[very thick, blue] (10,2) -- (10,4) node[pos = 0.5, right]{\small 3};

\draw[very thick, forest] (6,6) to[bend left = 45] (10,2);
\draw[very thick, forest] (12,6) to[bend right = 45] (8,2);
\draw[very thick, forest] (6,0) to[bend right = 45] (10,4);
\draw[very thick, forest] (12,0) to[bend left = 45] (8,4);
\draw[forest] (7.5,5.5) node{\small 4};
\draw[forest] (10.5,5.5) node{\small 4};
\draw[forest] (7.5,0.5) node{\small 4};
\draw[forest] (10.5,0.5) node{\small 4};

\filldraw (6, 0) circle (0.05 cm) node[left]{$b_3$};

\filldraw (6, 6) circle (0.05 cm)node[left]{$a_1$};

\filldraw (12, 0) circle (0.05 cm) node[right]{$a_4$};

\filldraw (12, 6) circle (0.05 cm) node[right]{$b_2$};

\filldraw (8, 2) circle (0.05 cm) node[below]{$a_3$};

\filldraw (8, 4) circle (0.05 cm) node[above]{$b_1$};

\filldraw (10, 2) circle (0.05 cm) node[below]{$b_4$};

\filldraw (10, 4) circle (0.05 cm) node[above]{$a_2$};

\end{tikzpicture}

\caption{The unique coloring of $H_4/2$, viewed as $K_{4,4}$ and as a folded hypercube}\label{folded 3 cube}

\end{center}

\end{figure}

\begin{thm}\label{thm: unique Pw+1-free coloring of Hw/2}
For $w\geq 3$, there is a unique (up to relabeling) rainbow $P_{w+1}$-free edge-coloring of $H_w/2$.   
\end{thm}
\begin{proof}
If $w \geq 5$, this follows quickly from Lemmas \ref{lem:H_w coloring exists},\ref{lem:H_w has TBFRC}, and \ref{lem:TBFRC extends to H_w}. Indeed, suppose $c$ and $c'$ are two rainbow $P_{w+1}$-free edge-colorings of $H_w/2$. (By Lemma~\ref{lem:H_w coloring exists}, at least one rainbow $P_{w+1}$-free edge-colorings of $H_w/2$ exists.) By Lemmas~\ref{lem:H_w has TBFRC} and \ref{lem:TBFRC extends to H_w}, each of $c$ and $c'$ contain a TBFRC, and in fact every TBFC in $H_w/2$ is rainbow under both $c$ and $c'$. Fix in particular the TBFC in $H_w/2$ whose $i$th edge corresponds to bit-flip $i$, which we shall call $C = (x_1, x_2, \dots, x_{w})$. Relabel the colors of $c'$ so that $c(x_ix_{i+1}) = c'(x_{i}x_{i+1})$ for every $i \in \{1,\dots, w\}$. By Lemma~\ref{lem:TBFRC extends to H_w}, every color which appears in $c'$ appears on some edge of $C$, so this relabeling of $c'$ indeed extends to a relabeling of every edge-color under $c'$. Now, under the described relabelling, $c$ and $c'$ agree on a TBFRC, so by Lemma~\ref{lem:TBFRC extends to H_w}, $c$ and the relabelling of $c'$ must be equal.

Thus, it remains only to verify the result for $w \in \{3,4\}$. If $w =3$, then we have $H_3/2 = K_4$. The only proper edge-coloring of $K_4$ that avoids a rainbow $P_4$-copy is the $3$-coloring under which each color class forms a perfect matching.

For $w = 4$, observe that $H_4/2 = K_{4,4}$. We show that there exists one rainbow $P_5$-free edge-coloring of $K_{4,4}$ up to isomorphism. Let $A \cup B$ be the bipartition of $V(K_{4,4})$ where every pair $a_i \in A, b_j \in B$ forms an edge $a_ib_j$. It is not difficult to check that, to avoid a rainbow $P_5$-copy, $K_{4,4}$ must be $4$ edge-colored. Thus, say we edge-color $K_{4,4}$ using colors $\{1,2,3,4\}$. Without loss of generality, let $c(a_1b_j) = j$ and $c(a_ib_1) = i$.

Consider the edge $a_ib_i$ for $i \geq 2$. If $c(a_ib_i) = j \neq 1$, then there exists a path $a_ib_ia_1b_1a_k$ for $k \neq j$ that is rainbow. Therefore $c(a_ib_i) = 1$.

Consider the edge $a_ib_j$ for $i,j \geq 2, i \neq j$. To maintain a proper edge-coloring, $c(a_ib_j)$ is not $i, j$, or $1$, and thus $c(a_ib_j)$ must be the unique color in $[4]
\setminus\{1,i,j\}$. Thus, up to relabelling, there is a unique rainbow $P_5$-free edge-coloring of $K_{4,4}$. For reference, we depict the described edge-coloring in Figure~\ref{folded 3 cube}.
\end{proof}

Now that we have established a precise understanding of the rainbow $P_{w+1}$-free colorings of $H_w/2$, we are ready to define and study a set of constructions arising from the folded hypercubes. Theorem~\ref{thm: unique Pw+1-free coloring of Hw/2} implies that, while it is possible to properly edge-color $H_w/2$ while avoiding a rainbow $P_{w+1}$-copy, any such edge-coloring contains many rainbow $P_w$-copies. In particular, every sequence of $w-1$ distinct bit-flips corresponds to a rainbow $P_w$-copy in $H_w/2$. In the following observation, we note that it is possible to construct such paths from any desired starting point, and which avoid any desired vertex or bit-flip.

\begin{obs}\label{obs: nice Pw in Hw/2}
    Let $w \geq 3$, $j \in [w]$, and $u,v \in V(H_w/2)$ be distinct vertices. Then there exist paths $P$ and $P'$ of length $w-1$ starting at $u$ in which each edge corresponds to a distinct bit-flip, such that no edge in $P$ corresponds to bit-flip $j$ and $v$ is not a vertex in $P'$. 
\end{obs}

Using this observation, we will show the desired upper bound on $\mathrm{sat}^*(n,P_k)$. Consider the following construction. 

\begin{const}
Fix $k \geq 6$ and $n \geq (k-1)2^{k-4}$. Let $G_k(n)$ be the graph obtained from $H_{k-3}/2$ by attaching $(k-2)$ pendant edges to each vertex of $H_{k-3}/2$ except for $0^{k-3}$, and $n - (k-1)2^{k-4} + k - 2$ pendant edges to $0^{k-3}$.
\end{const}

\begin{figure}[h!]
\begin{center}

\begin{tikzpicture}

\draw (0,0) -- (2,0) -- (2,2) -- (0,2) -- (0,0);
\draw (2,2) -- (0,0);
\draw (0,2) -- (2,0);

\draw (0,0) -- (-1.3,-1.5);
\draw (0,0) -- (-1.1,-1.7);
\draw (0,0) -- (-0.9,-1.9);
\draw (0,0) -- (-0.7,-2.1);

\draw (0,2) -- (-1.3,3.5);
\draw (0,2) -- (-1.1,3.7);
\draw (0,2) -- (-0.9,3.9);
\draw (0,2) -- (-0.7,4.1);

\draw (2,2) -- (3.3,3.5);
\draw (2,2) -- (3.1,3.7);
\draw (2,2) -- (2.9,3.9);
\draw (2,2) -- (2.7,4.1);

\draw (2,0) -- (3.3,-1.5);
\draw (2,0) -- (3.1,-1.7);
\draw (2,0) -- (2.9,-1.9);
\draw (2,0) -- (2.7,-2.1);

\filldraw (0,0) circle (0.05 cm);
\filldraw (2,0) circle (0.05 cm);
\filldraw (0,2) circle (0.05 cm);
\filldraw (2,2) circle (0.05 cm);

\filldraw (-1.3,-1.5) circle (0.05 cm);
\filldraw (-1.1,-1.7) circle (0.05 cm);
\filldraw (-0.9,-1.9) circle (0.05 cm);
\filldraw (-0.7,-2.1) circle (0.05 cm);

\filldraw (-1.3, 3.5) circle (0.05 cm);
\filldraw (-1.1, 3.7) circle (0.05 cm);
\filldraw (-0.9, 3.9) circle (0.05 cm);
\filldraw (-0.7, 4.1) circle (0.05 cm);

\filldraw (3.3,-1.5) circle (0.05 cm);
\filldraw (3.1,-1.7) circle (0.05 cm);
\filldraw (2.9,-1.9) circle (0.05 cm);
\filldraw (2.7,-2.1) circle (0.05 cm);

\filldraw (3.3, 3.5) circle (0.05 cm);
\filldraw (3.1, 3.7) circle (0.05 cm);
\filldraw (2.9, 3.9) circle (0.05 cm);
\filldraw (2.7, 4.1) circle (0.05 cm);

\draw (1, -3) node{(a)};

\draw (6,0) -- (8,0) -- (8,2) -- (6,2) -- (6,0);
\draw (8,2) -- (6,0);
\draw (6,2) -- (8,0);

\draw (6,0) -- (4.7,-1.5);
\draw (6,0) -- (4.9,-1.7);
\draw (6,0) -- (5.1,-1.9);
\draw (6,0) -- (5.3,-2.1);
\draw (6,0) -- (4.8,-1.6);
\draw (6,0) -- (5,-1.8);
\draw (6,0) -- (5.2,-2);
\draw (6,0) -- (5.4,-2.2);
\draw (6,0) -- (5.5,-2.3);
\draw (6,0) -- (5.65,-2.1);
\draw (6,0) -- (4.65,-1.4);

\draw (6,2) -- (4.7,3.5);
\draw (6,2) -- (4.9,3.7);
\draw (6,2) -- (5.1,3.9);
\draw (6,2) -- (5.3,4.1);

\draw (8,2) -- (9.3,3.5);
\draw (8,2) -- (9.1,3.7);
\draw (8,2) -- (8.9,3.9);
\draw (8,2) -- (8.7,4.1);

\draw (8,0) -- (9.3,-1.5);
\draw (8,0) -- (9.1,-1.7);
\draw (8,0) -- (8.9,-1.9);
\draw (8,0) -- (8.7,-2.1);

\filldraw (6,0) circle (0.05 cm);
\filldraw (8,0) circle (0.05 cm);
\filldraw (6,2) circle (0.05 cm);
\filldraw (8,2) circle (0.05 cm);

\filldraw[fill = blue!20!white] (5.2,-1.7) ellipse (0.6 and 0.75);

\draw (5.2, - 1.7) node{\small $n-16$};

\filldraw (4.7, 3.5) circle (0.05 cm);
\filldraw (4.9, 3.7) circle (0.05 cm);
\filldraw (5.1, 3.9) circle (0.05 cm);
\filldraw (5.3, 4.1) circle (0.05 cm);

\filldraw (9.3,-1.5) circle (0.05 cm);
\filldraw (9.1,-1.7) circle (0.05 cm);
\filldraw (8.9,-1.9) circle (0.05 cm);
\filldraw (8.7,-2.1) circle (0.05 cm);

\filldraw (9.3, 3.5) circle (0.05 cm);
\filldraw (9.1, 3.7) circle (0.05 cm);
\filldraw (8.9, 3.9) circle (0.05 cm);
\filldraw (8.7, 4.1) circle (0.05 cm);

\draw (7, -3) node{(b)};

\end{tikzpicture}

\caption{(a): The construction $G_6(n)$ for $n = 20$; (b): for general $n$ }\label{const 3 picture}
\end{center}
\end{figure}
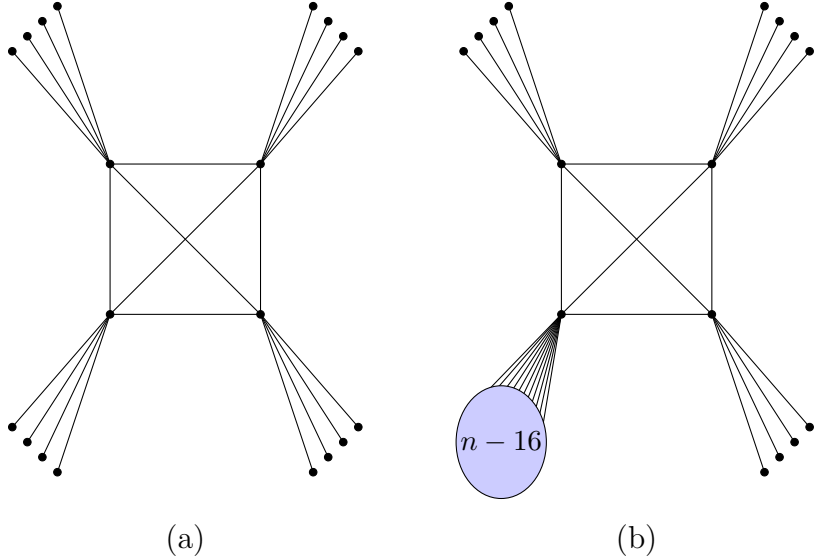

We call the subgraph of $G_k(n)$ isomorphic to $H_{k-3}/2$ the \textit{core} of $G_k(n)$, and call a vertex of $G_k(n)$ not in the core a \textit{pendant vertex}. We illustrate $G_6(20)$ and a more general depiction of $G_6(n)$ in Figure~\ref{const 3 picture}.

Note that when $n = (k-1)2^{k-4}$, every vertex of the core of $G_k(n)$ is incident to exactly $k-2$ pendant vertices. In general, we have
\begin{align*}
|E(G_k(n))| &= E(H_{k-3}/2) + (k-2)(2^{k-4} - 1) + n - (k-1)2^{k-4} + k - 2& \\
&= n + (k-3)2^{k-5} + (k-2)2^{k-4} -(k-2) - (k-1)2^{k-4} + k - 2& \\
&= n + (k-3)2^{k-5} - 2^{k-4}&\\
&= n + (k-5)2^{k-5}.&
\end{align*}

We will show that $G_k(n)$ is properly rainbow $P_k$-saturated in Theorem \ref{thm: sat*(Pk)<n+constant} with the help of the following simple lemma. 

\begin{lem}\label{lem: rb Pk in G iff rb Pk-2 in H}
    Let $k \geq 3$ let $C$ be an $n$-vertex graph with $n \geq 1$. Let $G$ be a graph obtained from $C$ by attaching at least $k-2$ pendant edges to each vertex of $C$. Then a proper edge-coloring $c$ of $G$ contains a rainbow $P_k$-copy if and only if the restriction of $c$ to $C$ contains a rainbow $P_{k-2}$-copy.
\end{lem}

\begin{proof}
    First, suppose $G$ contains a rainbow copy of $P_k$, say $P = (x_0, x_1, \dots, x_k)$.  Note that $P$ contains at most two pendant edges of $G$, and any pendant edge contained in $P$ must contain either $x_0$ or $x_k$. Thus, the sub-path $(x_1, \dots, x_{k-1})$ is a rainbow $P_{k-2}$-copy contained $C$.

    On the other hand, suppose $C$ contains a rainbow $P_{k-2}$-copy, say $P'$, with endpoints $u$ and $v$. Note that $u$ and $v$ are each incident to at least $k-2$ pendant edges in $G$. Furthermore, one of the $k-3$ edges on $P'$ is incident to $u$ and another is incident to $v$. So, only $k-4$ of the edge-colors used on $P'$ may appear among the $k-2$ pendant edges incident to $u$, and a (possibly different) set of $k-4$ of these colors may appear among the $k-2$ pendant edges incident to $v$. In particular, there are $2$ pendant edges incident to $u$ and two pendant edges incident to $v$ whose edge-colors are distinct from those on $P'$. We may extend $P'$ via two of these pendant edges to a rainbow $P_k$-copy in $G$.
\end{proof}

Observe that Lemma~\ref{lem: rb Pk in G iff rb Pk-2 in H} can be applied to $G_k(n)$, with the core of $G_k(n)$ acting as $C$; however, we will also apply Lemma~\ref{lem: rb Pk in G iff rb Pk-2 in H} in a case where $C$ may not equal the core of $G_k(n)$.

Now, we are ready to prove that $G_k(n)$ is properly rainbow $P_k$-saturated.

\begin{thm}\label{thm: sat*(Pk)<n+constant}
    For $k\geq 6$ and $n\geq (k-1)2^{k-4}$, $G_k(n)$ is rainbow $P_k$-saturated. In particular, $\sat^*(n,P_k) \leq n + (k-5)2^{k-5}$.
\end{thm}

\begin{proof}

    By Theorem \ref{thm: unique Pw+1-free coloring of Hw/2}, there is (up to relabelling of colors) a unique rainbow $P_{k-2}$-free coloring of the core, and by Lemma \ref{lem: rb Pk in G iff rb Pk-2 in H}, this extends to a rainbow $P_k$-free coloring of $G_k(n).$
    
    Next, we show that for any edge $e \not\in E(G_k(n))$,  any proper edge-coloring of $G_k(n) + e$ contains a rainbow $P_k$-copy. Assume, for the sake of contradiction, that some coloring $c$ of $G_k(n) + e$ is rainbow $P_k$-free for some $e \notin E(G_k(n))$. By Lemma \ref{lem: rb Pk in G iff rb Pk-2 in H}, the restriction of $c$ to the core must be rainbow $P_{k-2}$-free. By Theorem \ref{thm: unique Pw+1-free coloring of Hw/2}, the core must be colored so that edges corresponding to the same bit-flip receive the same color.

Now, either both endpoints of $e$ are pendant vertices, or else $e$ is incident to the core. Suppose first that both endpoints $x,y$ of $e$ are pendant vertices. Say $x,y$ have neighbors $u$ and $v$, respectively, in the core. As a consequence of Theorem~\ref{thm: unique Pw+1-free coloring of Hw/2}, the core is colored with precisely $k-3$ edge-colors, say $\{1,2,\dots, k-3\}$, and in particular, $u$ and $v$ are each incident to edges of every color $i \in \{1,2,\dots, k-3\}$ within the core. So $c(ux)$ and $c(vy)$ are not in $\{1,2,\dots, k-3\}$. Now, whether or not $c(xy) \in \{1,2, \dots, k-3\}$, Observation~\ref{obs: nice Pw in Hw/2} implies that there exists a rainbow $P_{k-3}$-copy $P$ in the core which has $u$ as an endpoint and does not contain an edge of color $c(xy)$. Let $w$ be the endpoint of this path not equal to $u$. Of at least $k-2$ pendant edges incident to $w$, at least $k-3$ are incident to neither $x$ nor $y$. As $k-3 \geq 3$, at least one such pendant edge, say $wz$, has $c(wz) \not \in \{c(xy), c(ux)\}$. Observe that the concatenation of $(y,x,u)$, $P$, and $(w,z)$ is a rainbow $P_k$-copy in $G_k(n) + e$, a contradiction.

Next, suppose $e$ is incident to the core. Let $C$ be the subgraph of $G_k(n)$ induced on the vertices of the core (so either $C$ is the core, or $C$ is the core with an added edge). By Lemma~\ref{lem: rb Pk in G iff rb Pk-2 in H}, to avoid a rainbow $P_k$-copy, we must color so that $C$ is rainbow $P_{k-2}$-free. We claim that this is not possible if both endpoints of $e$ are incident to the core. 

Indeed, to avoid an immediate contradiction, all core edges which are not equal to $e$ still must be colored as described by Theorem~\ref{thm: unique Pw+1-free coloring of Hw/2}. If both endpoints $u,v$ of $e$ are contained in the core, then by Observation~\ref{obs: nice Pw in Hw/2}, there exists a rainbow $P_{k-3}$-copy $P$ in the core which has $u$ as an endpoint and does not contain $v$. Note also that, since $u,v$ are already incident to edges of all colors from $\{1,2,\dots, k-3\}$, we must have $c(uv) \not\in \{1,2,\dots, k-3\}$. Thus, concatenating $uv$ with $P$ yields a rainbow $P_{k-2}$-copy in the core, which can be extended to a rainbow $P_k$-copy in $G_k(n) + e$ by Lemma~\ref{lem: rb Pk in G iff rb Pk-2 in H}, a contradiction. 

Finally, we consider the case where one endpoint $u$ of $e$ is in the core and the other, say $x$, is a pendant vertex. Note again that $c(ux) \not \in \{1,2, \dots, k-3\}$ and since $ux \not\in E(G_k(n)),$ $x$ is adjacent to a core vertex distinct from $u$, say $v$. Now, using Observation~\ref{obs: nice Pw in Hw/2}, there is a rainbow $P_{k-3}$-copy $P$ in the core with $v$ as an endpoint which does not contain $u$. Concatenating $P$ with $(u,x,v)$ and a pendant edge incident to the other endpoint of $P$ will again yield a rainbow $P_k$-copy in $G_k(n) + e$, a contradiction.
\end{proof}

Finally, for the sake of completeness, we derive Theorem~\ref{paths thm}

\pathsthm*

\begin{proof}
We have the lower bound by Proposition~\ref{lem:P5-lowerbound-acycliccomps} and the upper bound for $k\geq 6$ by Theorem~\ref{thm: sat*(Pk)<n+constant}. For $k = 5$, consider the $n$-vertex graph $G_5(n)$ obtained by attaching $n-4$ pendant edges to one vertex of a $K_4$-copy; we again call this $K_4$-copy the \textit{core} of $G_5(n)$ and say that vertices outside the core are \textit{pendant vertices}. Note that we can extend the perfect matching coloring of the core to a rainbow $P_5$-free proper edge-coloring of $G_5(n)$. We claim that in fact, $G_5(n)$ is properly rainbow $P_5$-saturated.

Let $x$ be the core vertex to which all pendant vertices of $G_5(n)$ are adjacent. Since $n \geq 4\cdot 2  = 8$, $x$ has at least $4$ pendant vertex neighbors. Since the core is complete, any edge added to $G_5(n)$ either connects two pendant vertices $u,v$ or connects a core vertex $y \neq x$ to a pendant vertex $u$. Note that if we add edge $uv$ to $G_5(n)$, then $\{u,v,x\}$ must span a (rainbow) triangle, while if we add edge $uy$ to $G_5(n)$, then there exist two other pendant vertices $v,w$ such that $c(ux), c(vx), c(wx),$ and $c(uy)$ are pairwise distinct. We depict these two cases in Figure~\ref{k=5 cases}, also labeling the other core vertices of $G_5(n)$.

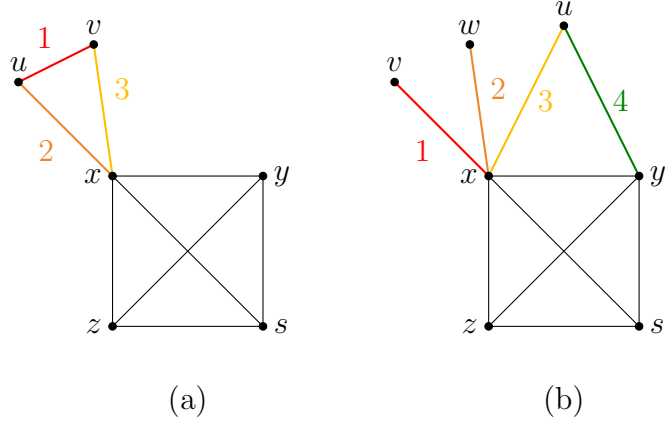
\begin{figure}[h!]
\begin{center}

\begin{tikzpicture}

\draw (0,0) -- (2,0) -- (2,2) -- (0,2) -- (0,0);
\draw (0,0) -- (2,2);
\draw (2,0) -- (0,2);

\draw[thick, red] (-1.25,3.25) -- (-0.25, 3.75) node[pos = 0.6, above left]{1};
\draw[thick, cadmium] (-1.25,3.25) -- (0,2) node[pos = 0.5, below left]{2};
\draw[thick, amber] (-0.25, 3.75) -- (0,2) node[pos = 0.5, above right]{3};

\filldraw (0,0) circle (0.05 cm) node[left]{$z$};
\filldraw (2,0) circle (0.05 cm) node[right]{$s$};
\filldraw (0,2) circle (0.05 cm) node[left]{$x$};
\filldraw (2,2) circle (0.05 cm) node[right]{$y$};
\filldraw (-1.25,3.25) circle (0.05 cm) node[above]{$u$};
\filldraw (-0.25,3.75) circle (0.05 cm) node[above]{$v$};

\draw (1,-1)  node{(a)};

\draw (5,0) -- (7,0) -- (7,2) -- (5,2) -- (5,0);
\draw (5,0) -- (7,2);
\draw (7,0) -- (5,2);

\draw[thick, red] (3.75,3.25) -- (5,2) node[pos = 0.5, below left]{1};
\draw[thick, cadmium] (4.75, 3.75) -- (5,2) node[pos = 0.5, above right]{2};

\draw[thick, amber] (5,2) -- (6,4) node[pos = 0.5, right]{3};

\draw[thick, forest] (6,4) -- (7,2) node[pos = 0.5, right]{4};

\filldraw (5,0) circle (0.05 cm) node[left]{$z$};
\filldraw (7,0) circle (0.05 cm) node[right]{$s$};
\filldraw (5,2) circle (0.05 cm) node[left]{$x$};
\filldraw (7,2) circle (0.05 cm) node[right]{$y$};
\filldraw (3.75,3.25) circle (0.05 cm) node[above]{$v$};
\filldraw (4.75,3.75) circle (0.05 cm) node[above]{$w$};
\filldraw (6,4) circle (0.05 cm) node[above]{$u$};

\draw (6,-1)  node{(b)};

\end{tikzpicture}

\caption{(a): $G_5(n)$ with $uv$ added between pendant vertices; (b): $G_5(n)$ with $uy$ added between a pendant vertex and a core vertex}\label{k=5 cases}

\end{center}
\end{figure}

In $G_5(n) + uv$, note that at most one of $c(xy), c(xz), c(xs)$ is equal to $1$; without loss of generality, $c(xy) \neq 1$. Thus, both $(u,v,x,y)$ and $(v,u,x,y)$ are rainbow $P_4$-copies in $G_5(n) + uv$. Note also that $c(yz), c(ys)$ are not equal to $c(xy)$, and at most one is equal to $1$. Without loss of generality, $c(yz) \neq 1$. Then either $(u,v,x,y,z)$ or $(v,u,x,y,z)$ is a rainbow $P_5$-copy in $G_5(n) + uv$. Thus, any proper edge-coloring of $G_5(n) + uv$ contains a rainbow $P_5$-copy.

In $G_5(n) + uy$, note that at most one of $c(yz), c(ys)$ is equal to $3$. Without loss of generality, $c(yz) \neq 3$. Then either $(v, x, u, y, z)$ or $(w, x, u, y, z)$ is a rainbow $P_5$-copy in $G_5(n) + uy$. Thus, any proper edge-coloring of $G_5(n) + uy$ contains a rainbow $P_5$-copy.

We conclude that $G_5(n)$ is properly rainbow $P_5$-saturated. Note that $|E(G_5(n))| = n + 2$. Thus, the desired upper bound also holds for $k = 5$.
\end{proof}

\section{Cycles}\label{sec:cycles}

\cyclesthm*

\begin{proof}
    Fix odd $k \geq 7$, and suppose $n \geq 3k - 2$.     
    Let $G = X \vee Y$, where $X = K_{\frac{k-1}{2}}$ and $Y = E_{n-\frac{k-1}{2}}$. We will show that $G$ is rainbow $C_k$-saturated. 
    
    Set $V(X) = \{x_1, \dots, x_{\frac{k - 1}{2}}\}$ and $V(Y) = \{y_1, \dots, y_{n - |V(X)|}\}$.  Since $|V(X)| = \frac{k - 1}{2}$ and $Y$ is an independent set, there are no odd cycles of length $k$ or greater in $G$. Let $G' = G + y_1y_2$ and $c: E(G') \to \mathbb{N}$ be a proper edge-coloring of $G'$. Observe that the double star $S(1, 2)$ is a subgraph of $G'[x_1, x_2, x_3, y_1, y_2]$ with centers $y_1, y_2$. So, there is a rainbow path of length $3$ containing $y_1y_2$. Without loss of generality, suppose this rainbow path is $P = x_1y_1y_2x_2$. \\
    
    We begin building our rainbow $C_k$ by setting $C := P$ and greedily building the rest. For an edge $e \in G'$, we call $e$ a bad edge if $c(e) \in c(C)$ and otherwise, $e$ is a good edge. Set
    \begin{equation*}
        S_3 := \{y \in Y : y \not \in V(C) \text{ and } x_2y, x_3y \text{ are good edges}\}.
    \end{equation*}
    
    Observe that there can be at most two bad edges outside of $E(C)$ incident to $x_2$ and at most three bad edges incident to $x_3$. Therefore, if $|V(Y)| \geq 2 + 2 + 3 + 1 = 8$, then $S_3 \neq \emptyset$. Since $|V(Y)| \geq \frac{5}{2}(k - 1) + 1 \geq 16$, this inequality holds. By relabeling, we may assume $y_3 \in S_3$. Now, add the edges $x_2y_3, y_3x_3$ to $C$. Continuing in this way, after relabeling, suppose $C = x_1y_1y_2x_2y_3x_3 \dots y_ix_i$, $i < \frac{k - 1}{2}$. Then, set
    \begin{equation*}
        S_{i + 1} := \{y \in V(Y) : y \not \in V(C) \text{ and } x_iy, x_{i+1}y \text{ are good edges}\}.
    \end{equation*}

    As before, observe that there can be at most $|E(C)| - 1 = 2i$ bad edges outside of $E(C)$ incident to $x_i$ and at most $|E(C)| = 2i + 1$ bad edges incident to $x_{i+1}$. Therefore, $S_{i + 1} \neq \emptyset$ if the following inequality holds:
    \begin{equation*}
        |V(Y)| \geq |V(C) \cap V(Y)| + (2i) + (2i + 1) + 1 = i + (2i) + (2i + 1) + 1 = 5i + 2
    \end{equation*}

    Since $|V(Y)|\geq \frac{5}{2}(k - 1) + 1$, this inequality holds for all $3 \leq i < \frac{k - 1}{2}$ and we may continue building $C$. Finally, suppose $C = x_1y_1y_2x_2y_3x_3 \dots y_{\frac{k - 1}{2}}x_{\frac{k - 1}{2}}$ and set
    \begin{equation*}
        S_{\frac{k - 1}{2}} := \{y \in V(Y) : y \not \in V(C) \text{ and } x_1y, x_{\frac{k - 1}{2}}y \text{ are good edges}\}.
    \end{equation*}

    Now, there can be at most $|E(C)| - 1 = k - 1$ bad edges outside of $E(C)$ incident to $x_1$ and same for $x_{\frac{k - 1}{2}}$. Therefore, $S_{\frac{k - 1}{2}} \neq \emptyset$ if the following inequality holds:
    \begin{equation*}
        |V(Y)| \geq |V(C) \cap V(Y)| + (k - 1) + (k - 1) + 1 = \frac{k - 1}{2} + 2(k - 1) + 1 = \frac{5}{2}(k - 1) + 1.
    \end{equation*}

    Again, since $|V(Y)| \geq \frac{5}{2}(k - 1) + 1$, the inequality holds and we may find $y_{\frac{k + 1}{2}} \in S_{\frac{k - 1}{2}}$, add the edges $x_1y_{\frac{k + 1}{2}}$ and $x_{\frac{k - 1}{2}}y_{\frac{k + 1}{2}}$ to $C$, completing our rainbow $C_k$. Therefore, $G'$ is rainbow $C_k$-saturated. 

    Finally, we count the edges of $G'$ as follows:
    \begin{equation*}
        \begin{split}
            |E(G')|
                &=  |E(X)| + |E(X, Y)| + |E(Y)| \\
                &=  \binom{\frac{k - 1}{2}}{2} + |V(X)|(n - |V(X)|) + 0 \\
                &=  \binom{\frac{k - 1}{2}}{2} + \left( \frac{k -1}{2} \right) \left(n - \frac{k - 1}{2} \right) \\
                & = \left( \frac{k - 1}{2} \right) n -\binom{\frac{k + 1}{2}}{2}. \\ 
        \end{split}
    \end{equation*}
\end{proof}

\section*{Acknowledgements}
Work on this project started during the Research Training Group (RTG) rotation at Iowa State University in the spring of 2024. Dustin Baker, Enrique Gomez-Leos, Emily Heath, Joe Miller, Hope Pungello, and Nick Veldt were supported by NSF grant DMS-1839918.  Ryan Martin was supported by Simons Collaboration Grant \#709641. 

\bibliography{references}
\bibliographystyle{abbrv}

\end{document}